\documentclass[11pt]{article}

\usepackage{graphicx}
\usepackage{amsmath}
\usepackage{amssymb}
\usepackage{amsthm}
\usepackage{hyperref}
\usepackage[margin=1in]{geometry}
\usepackage{xcolor}
\usepackage{caption}
\usepackage{subcaption}
\usepackage{lineno}
\usepackage[numbers, sort&compress, square]{natbib}

\newcommand{\email}[1]{\href{mailto:#1}{#1}}
\newtheorem{remark}{Remark}
\newtheorem{theorem}{Theorem}
\newtheorem{lemma}[theorem]{Lemma}

\usepackage{multirow}

\title{A Convergent Geometry-Aware Reduction for Diffusion in Branched Tubular Networks}
\author{Zachary M. Miksis  \thanks{Department of Mathematics, Temple University, Philadelphia, Pennsylvania, USA  (\email{miksis@temple.edu}, \email{queisser@temple.edu})} \and Gillian Queisser\footnotemark[1]}
\date{}

\begin{document}

\maketitle

\begin{abstract}
    Diffusion in tubular networks with variable radius arises in biological, engineering, and physical applications. The standard one-dimensional reduction is the Fick-Jacobs equation, which Zwanzig showed in 1992 is unstable and inaccurate when the radial gradient is large; three decades of subsequent modifications have not resolved this. We show by elementary analysis that these modifications inherit a structural inconsistency from truncation in Jacobs's original derivation, and we re-derive the model from the three-dimensional axisymmetric diffusion equation by Taylor expansion at finite subdomain width. The result is a geometry-aware expansion of the Fick-Jacobs equation with model error independent of the radial gradient, a finite-volume discretization with provable discrete energy stability under a small-grid threshold and provable convergence at the rate $\mathcal{O}(\Delta x^2)$, and natural treatment of branched networks at equivalent computational cost. On the canonical truncated-cone benchmark of Berezhkovskii et al., a tunable form of the expansion model reproduces axisymmetric ground-truth mean first-passage times across the full slope range, including the regime $|R_x| > 1$ where no parameter choice in prior corrections can achieve the same agreement. Applied to a dendritic spine, the one-dimensional reduction recovers the qualitative calcium dynamics of a full three-dimensional simulation at a fraction of the computational cost.
\end{abstract}



\section{Introduction}

 Many biological, engineering, and experimental applications fundamentally rely on diffusion processes in bounded domains, particularly with porous boundaries. Among many others, such processes arise in ion transport through biological membranes \cite{Hille2001,YKA1991}, carbon nanotubes \cite{BH2002}, diffusion across a porous medium interface \cite{ARDIR2014}, biosensors \cite{Keyser2006}, or permeation tubes \cite{Lucero1971}. For any application in which the domain geometry can affect the underlying (bio)physical processes it becomes necessary to incorporate the geometric features into the mathematical model. A straightforward way is to cast and solve the problem in the full three-dimensional domain. 
 The examples mentioned above describe processes in channels, tubes, and tubular networks and many methods are available to computationally treat these types of problems accurately, including Monte Carlo simulations \cite{SLK2012} and finite element based methods \cite{GQ2022}. However, computational bottlenecks arise when discretizing large tubular network structures that introduce a very large number of degrees of freedom, as would happen, e.g., in structures present in the cardiovascular system or the brain. In three space dimensions, this becomes computationally expensive and simulations quickly become intractable. 

To address the problem posed by computational cost, methods based on the Fick-Jacobs equation \cite{Jacobs1935} are frequently proposed. This equation reduces the diffusion equation in a three-dimensional tube to a one-dimensional equation that respects the variable radius of the domain, opening the door for spatially larger or temporally longer simulations. It was briefly derived by Jacobs in a single paragraph of his classical textbook, without expansion of the full mathematical details. His approach followed standard continuum derivation methodology used to derive Fourier's heat equation and Fick's diffusion law. While sound for uniform domains, it loses essential geometric information when applied to domains of varying radius, as we show. 

Zwanzig later showed through rigorous analysis \cite{Zwanzig1992} that this model diverges from the exact solution in a hyperboloidal cone in the regime $|R_x| \gtrsim 1$. Notably, the Fick-Jacobs equation is motivated precisely by domains of significant geometric variation, the regime where this instability is most pronounced. In domains of weak geometric variation, approximation by the standard diffusion equation may itself be sufficient. To extend the region of stability, he proposed a spatially adaptive diffusion coefficient based on the changing radius. Other adaptive diffusion coefficients have been proposed since then \cite{RR01, KP06}, as well as temporal corrections \cite{Kalinay2013} and the inclusion of biharmonic coefficients \cite{CLG2023}. While they improve on the stability of the classic Fick-Jacobs equation, it has been noted that many still present computational instabilities due to the radial gradient \cite{BPB2007}. As we demonstrate in this work, these corrections all inherit a common structural inconsistency from truncation in Jacobs' original derivation that the standard diffusion coefficient modifications do not address.

 An additional complication arises in branching networks. While the previously mentioned methods are applicable to long tubes, they do not computationally address how to treat branching points in tubes that form large networks. A common method for addressing this is via the Hines matrix \cite{Hines1984} which is based on an indexing strategy of the network nodes. However, this requires an additional reordering of the network done as a preconditioning step. A more direct method is to discretize the branch nodes in a way that respects conservation principles, as was done in \cite{BRNQ2023} for diffusion problems.
 
A geometry-aware model retains 3D geomoetric information, such as radius and radial gradient, in a 1D framework. In this paper, we return to the foundations of Jacobs' derivation \cite{Jacobs1935} and treat it as a Taylor expansion, identifying two sources of geometric information loss that are inherited by all previous standard corrections. We address these directly through an expansion of the computed flux, producing a geometry-aware model with a provably stable discretization which converges under refinement to the geometry-aware reduction it is derived from, whereas the standard diffusion-coefficient corrections, by an elementary argument we give in Sec.~\ref{sec:original-fick-jacobs}, cannot. We note that both the expanded model and the classical Fick-Jacobs equation share the same continuous limit, with the distinction lying entirely in geometric fidelity at the level of the finite subdomain, where all computation occurs. We also provide simple and computationally efficient discretizations for branching nodes, and demonstrate the accuracy and stability of our method in numerical experiments.

Achieving computational efficiency in these problems has also been pursued through neural network approaches \cite{GG2024, MQ2024, MFK2021}, though these methods face limitations in flexibility and applicability to branched networks. The dimension-reduced model developed here achieves efficiency through mathematically rigorous model reduction rather than learned approximation, preserving geometric fidelity while enabling simulations that would otherwise require full three-dimensional computation. We demonstrate how this can be particularly impactful in biological applications where tubular network structures are ubiquitous and long-timescale simulations are essential for clinical relevance, opening the door to parameter studies, clinical timescale simulations, and multiscale modeling that remain computationally intractable in three dimensions.

The structure of the paper is as follows. In Section 2, we present the classical models of 1D diffusion and discuss their drawbacks, including an elementary analysis of Jacobs' original derivation that identifies the structural sources of error inherited by all standard corrections. In Section 3, we derive our expansion of the Fick-Jacobs model as a systematic reduction of the 3D axisymmetric diffusion equation on the truncated cone, in which the classical Fick-Jacobs equation arises as the leading-order term and our model retains the next order in the subdomain width $\Delta x$, and incorporate lateral fluxes in a natural and simplified manner. In Section 4, we describe the finite volume discretization of our model and the application of branch nodes in a one-dimensional network. In Section 5, we establish discrete energy stability and convergence of the numerical method. Section 6 presents numerical experiments, including a manufactured-solution convergence test, a benchmark against axisymmetric ground truth on the truncated cone of Berezhkovskii et al. \cite{BPB2007}, and a biological application to a dendritic spine.



\section{Classical Models}

\begin{figure}
    \centering
    \begin{subfigure}[b]{0.48\textwidth}
        \includegraphics[scale=0.7]{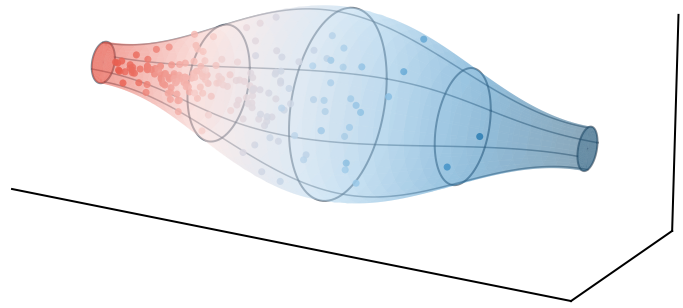}
        \caption{}
    \end{subfigure}
    \begin{subfigure}[b]{0.48\textwidth}
        \includegraphics[scale=0.7]{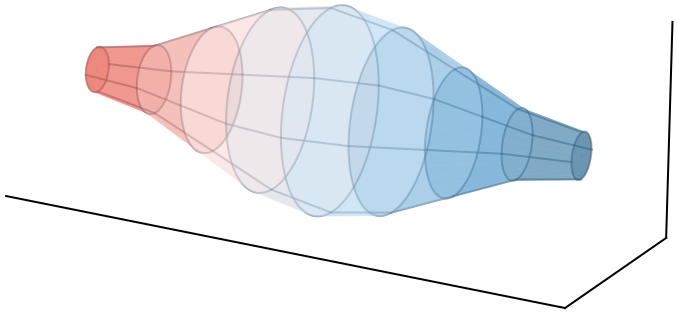}
        \caption{}
    \end{subfigure}
    \caption{Diffusion through a 3D channel domain of varying radius. (a) Smooth 3D domain with diffusive particles. (b) 3D discretized domain divided into truncated cone subdomains.}
    \label{fig:standard}
\end{figure}

In this section, we present the classical models of 1D diffusion, as well as lay out their drawbacks and model corrections from the literature. Each approach aims to construct a reduced 1D model of diffusion through a 3D channel of varying radius, such as in the channel shown in Fig.~\ref{fig:standard}. The general process is to divide the larger domain into smaller subdomains, such as those in the middle region of the figure, and solve over each subdomain. 

\subsection{1D Standard Reduction}

The standard 1D reduction of diffusion through a tube, used in an array of applications such as cable theory \cite{BG2021} and ionic transport \cite{BRNQ2023}, is the classical Fick's law,
\begin{equation}
    \frac{\partial u}{\partial t} = D_0 \frac{\partial^2 u}{\partial x^2}. \label{eq:fick}
\end{equation}
This treats each subdomain as a cylinder of constant radius, where the concentration of interest $u$ is assumed to be radially uniform. Discretizing this effectively means treating each subsection (volume) as its own cylinder, resulting in a ``stack'' of cylinders. This works well when averaging over domains where the length scale is much larger than the radius, and the radial derivative is near zero. Where this can fail is when the radial derivative is large. What results is stacked cylinders that are mismatched in radius. This reduction then loses or misses dynamics that may be attributed to the geometric structure of the domain, establishing the key problem we address here.


\subsection{Fick-Jacobs Model}\label{sec:original-fick-jacobs}
\begin{figure}[h!]
    \centering
    \includegraphics[scale=0.2]{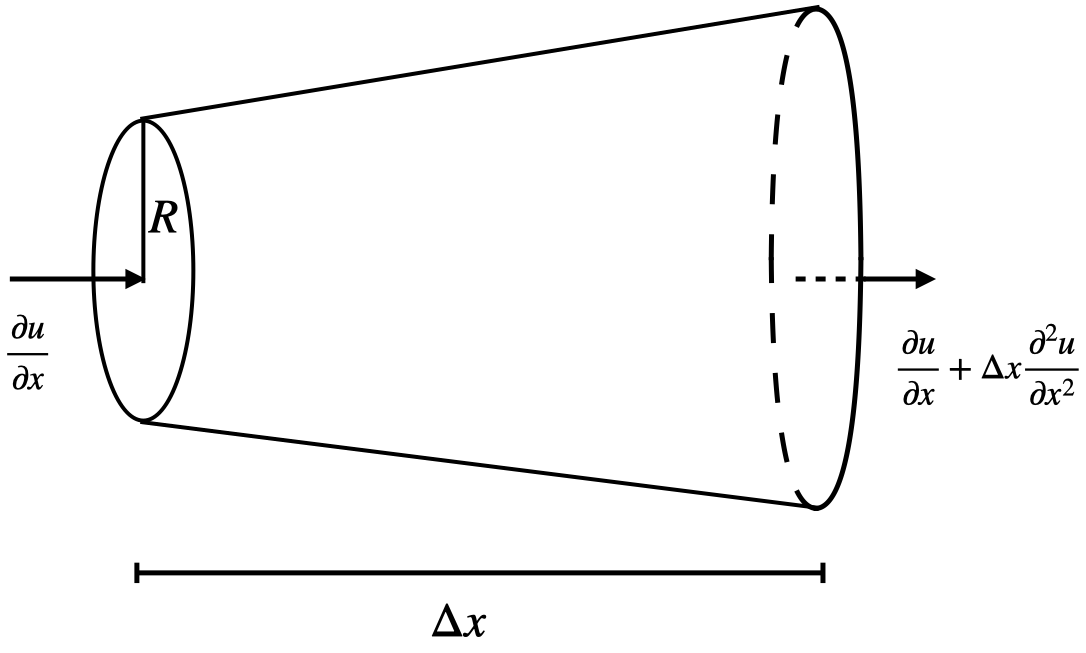}
    \caption{Jacobs flux through a truncated cone.\textbf{}}
    \label{fig:jacobs-flux}
\end{figure}

In \cite{Jacobs1935}, Jacobs provided a derivation of a 1D reduced model of diffusion through a domain with varying radius, specifically of the truncated cone in Fig.~\ref{fig:jacobs-flux}. This corresponds to a single subdomain, with the basic assumption of radial and rotational concentration symmetry in an infinitesimal domain.
As the domain is exactly a truncated cone, it can be assumed that $R_x$ is a constant value and $R_{xx} = 0$. It is also assumed that the change in volumetric concentration over time is equal to the difference between the entering and exiting fluxes.

Taking the total flux through the left face of the truncated cone and subtracting the total flux through the right face, we get the total flux through the domain. Setting this equal to the volumetric flux and dividing both sides by $\pi \Delta x R^2$,
\begin{equation}
    \left( 1 + \Delta x \frac{R_x}{R} + \frac{\Delta x^2}{3} \frac{R_x^2}{R^2} \right) \frac{\partial u}{\partial t} = \frac{D_0}{R^2} \frac{\partial}{\partial x} \left( R^2 \frac{\partial u}{\partial x} \right) + D_0 \Delta x \left( \frac{R_x^2}{R^2} \frac{\partial u}{\partial x} + 2 \frac{R_x}{R} \frac{\partial^2 u}{\partial x^2} \right) + D_0 \Delta x^2 \frac{R_x^2}{R^2} \frac{\partial^2 u}{\partial x^2}. \label{eq:fick-jacobs-exp}
\end{equation}
 This complete form of the equation never appeared in Jacobs' original work, in particular the volumetric term on the left hand side. At this point, Jacobs drops infinitesimal $\Delta x$ terms and derived what is now referred to as the Fick-Jacobs equation,
\begin{equation}
    \frac{\partial u}{\partial t} = \frac{D_0}{R^2} \frac{\partial}{\partial x} \left( R^2 \frac{\partial u}{\partial x} \right), \label{eq:fick-jacobs}
\end{equation}
where $u$ denotes the concentration of some quantity of interest, $D_0$ the constant diffusion coefficient of $u$, and $R$ the domain's cross-sectional radius, which depends on the axial direction $x$. 

It can be noticed that Eq.~\eqref{eq:fick-jacobs-exp} is essentially a Taylor expansion around the left side of the domain, with right hand side error term
\begin{equation}
    \mathcal{E} = D_0 \Delta x \left( \frac{R_x^2}{R^2} \frac{\partial u}{\partial x} + 2 \frac{R_x}{R} \frac{\partial^2 u}{\partial x^2} \right) + D_0 \Delta x^2 \frac{R_x^2}{R^2} \frac{\partial^2 u}{\partial x^2}.
\end{equation}
This error term is dominated by and grows like $R_x^2$ when $|R_x| > 1$. This is a fundamental problem with the classical Fick-Jacobs model that modern modifications attempt to address.

Additionally, it may be noticed that the volumetric term on the left hand side is dropped entirely. In practice, this means that the left hand side is now effectively treated as the flux through a perfect cylinder, meaning that the advective term on the right hand side is no longer acting as a corrective term in radially varying geometry but as an additional source of error on a diffusion equation. This structural error is inherited by all subsequent modifications, as they are derived from the truncated equation. The diffusion coefficient modifications considered here do not address the loss of volumetric geometric information on the left hand side, as we demonstrate numerically in later sections.

We note that these additional terms arise from applying standard continuum derivation methodology to a non-uniform domain. When the domain is uniform, this approach is sound and dropping the higher order terms is perfectly valid. In non-uniform domains, these higher order terms contain essential geometric information that is lost in the truncation.


\subsubsection{Diffusion Coefficient Correction}

Jacobs derived his entire model in a single paragraph with two equations, one of which was the final result, which resulted in the error term going unnoticed. Zwanzig identified the model error through rigorous analysis \cite{Zwanzig1992}, but did not identify the fundamental source, namely the truncation of higher order terms. He demonstrated that the Fick-Jacobs equation diverged from the exact solution of diffusion in an infinitely long hyperboloidal cone by identifying that the radial equilibrium assumption fails in regimes of large radial gradient, but this was guaranteed to fail since Jacobs' derivation was in an infinitesimally small domain that removed volumetric information. He proposed a diffusion coefficient correction to account for the radial change that tracked the exact solution over a wider range of $R_x$ than the standard Fick-Jacobs model, but still ultimately diverged due to the missing geometric information. Zwanzig's work then guided the subsequent development of the field, as the community focused on developing diffusion coefficient corrections to address the error. Different models have been proposed to modify the diffusion coefficient and improve results on steep domains \cite{BDB2015}. Common ones include Zwanzig (Zw) \cite{Zwanzig1992}, Reguera-Rubi (RR) \cite{RR01}, Kalinay-Percus (KP) \cite{KP06}, and Dorfman-Yariv (DY) \cite{DY2014}:
\begin{eqnarray}
    D_{FJ}(x) &=& D_0 \\
    D_{Zw}(x) &=& \frac{D_0}{1 + [R_x]^2/2} \\
    D_{RR}(x) &=& \frac{D_0}{(1 + [R_x]^2)^{1/2}} \\
    D_{KP}(x) &=& \frac{D_0 \cdot \arctan([R_x]/2)}{[R_x]/2} \\
    D_{DY}(x) &=& D_0 \left( \frac{1}{1/R} - \frac{1}{3}\frac{[R_x]^2/R}{(1/R)^2} + \frac{4([R_x]^4/R}{45(1/R)^2} + \frac{([R_x]^2/R)^2}{9(1/R)^3} + \frac{R[R_{xx}]^2}{45(1/R)^2} \right)
\end{eqnarray}
Unlike our model presented in Sec.~\ref{sec:expanded-fick-jacobs}, these are mainly developed by different methods of deriving and analyzing the Fick-Jacobs equation from the Smoluchowski equation. It should be noted that Martens et al. \cite{Martens2011} arrived at the same correction as Kaylinay-Percus using a slightly simpler analysis, but still utilizing the Smoluchowski framework. A key aspect missing from these derivations is that the error in the original model is still dominated by $R_x$. We demonstrate this numerically in later sections.


\subsubsection{Kalinay Temporal Correction}\label{sec:kalinay}

Kalinay extended the Zwanzig derivation to add temporal correction terms similar to those we derive for our model in eq.~\eqref{eq:dFJ}, and demonstrated the improved effectiveness of this on a 2D cone \cite{Kalinay2013}. Following the process described in the next section for 2D problems results in the original Fick-Jacobs equation without higher order terms, so we will not address this here.

The Kalinay correction affects only the left side of the equation. It is done by setting the left hand side of the Fick-Jacobs equation to
\begin{equation}
    [1 + \partial_x g(x)]\frac{\partial u}{\partial t},
\end{equation}
where, letting $\epsilon$ be some scaling term,
\begin{equation}
    g(x) = \frac{x}{2} \left( \frac{\arctan \sqrt{\epsilon} R_x}{\sqrt{\epsilon} R_x} + \frac{\sqrt{\epsilon} R_x}{3} \arctan \sqrt{\epsilon} R_x - 1 \right).
\end{equation}
It is assumed that $\epsilon = 1$ in all numerical tests, as in \cite{Kalinay2013}.

While this correction modifies the left hand side of the equation, it does not recover the volumetric geometric information lost in Jacobs' original truncation. The correction to $g(x)$ arises from a temporal mapping in the Smoluchowski framework rather than from geometric consideration of the local volume, and therefore does not address the structural loss identified in Section~\ref{sec:original-fick-jacobs}.



\section{Expanded Fick-Jacobs Model}
\label{sec:expanded-fick-jacobs}

In this section, we derive our expanded form of the Fick-Jacobs equation. The derivation begins from the three-dimensional axisymmetric diffusion equation and proceeds by cross-sectional averaging followed by Taylor expansion of the resulting conservation law at finite subdomain width $\Delta x$, refining Jacobs' original flux-balance derivation \cite{Jacobs1935} by retaining higher-order terms in the expansion rather than truncating at first order. This is structurally distinct from the established slenderness-asymptotic family of reductions \cite{Zwanzig1992, RR01, KP06, DY2014}, which expands the 3D solution in a slenderness parameter $\epsilon = R/L$ and produces, at next order beyond Fick-Jacobs, a slope-dependent rescaling of the effective diffusivity, $D_{\mathrm{eff}}(R') = D_0[1 - (R')^2/2 + \mathcal{O}((R')^4)]$ in the Zwanzig form. The two reductions agree at leading order but differ at next order: the slenderness family produces only a rescaling of the leading-order flux, while the expansion model also produces a third-derivative flux $b(x)\partial_x^3 u$ and finite-$\Delta x$ coefficients on both the flux and the time derivative.

\label{sec:flux-expansion-derivation}

The starting point is the three-dimensional axisymmetric diffusion equation
\begin{equation}
    \partial_t u \;=\; D_0 \left[ \frac{1}{r}\partial_r(r \partial_r u) + \partial_x^2 u \right] \qquad \text{on } \{(x, r): 0 \leq r \leq R(x)\}, \label{eq:3d-axisym}
\end{equation}
posed on the truncated-cone subdomain $\Omega_{\Delta x} = \{(x, r, \theta): x_0 - \Delta x/2 \leq x \leq x_0 + \Delta x/2,\ 0 \leq r \leq R(x),\ 0 \leq \theta < 2\pi\}$ with lateral wall $r = R(x)$ on which the no-flux boundary condition $\hat n \cdot \nabla u = 0$ holds, where $\hat n = (-R_x, 1)/\sqrt{1 + R_x^2}$ is the outward normal in the $(x, r)$ plane. Following Jacobs \cite{Jacobs1935}, the subdomain is a truncated cone so $R(x) = R(x_0) + R_x(x_0)(x - x_0)$ is linear and $R_{xx} = 0$ on the subdomain.

Define the cross-sectional average $\bar u(x, t) := (2/R(x)^2)\int_0^{R(x)} r\, u(x, r, t)\, dr$. Multiplying \eqref{eq:3d-axisym} by $2\pi r$ and integrating over the cross-section $r \in [0, R(x)]$, the radial-flux divergence integrates exactly via the lateral boundary condition, producing a reduced identity for $\bar u$. Integrating once more over the subdomain $x \in [x_0 - \Delta x/2, x_0 + \Delta x/2]$ and applying the divergence theorem,
\begin{equation}
    \frac{d}{dt}\int_{x_0 - \Delta x/2}^{x_0 + \Delta x/2} \pi R(x)^2 \bar u(x, t)\, dx \;=\; \Phi(x_0 - \Delta x/2) - \Phi(x_0 + \Delta x/2) + \Phi_{\mathrm{lat}}, \label{eq:slab-conservation}
\end{equation}
where $\Phi(x) := -D_0 \pi R(x)^2 \partial_x \langle u \rangle(x, t)$ is the axial cap flux ($\langle u \rangle$ being the cross-sectional mass-weighted average of the 3D field $u$) and $\Phi_{\mathrm{lat}}$ is the lateral surface flux, which vanishes for the closed-wall problem but is reintroduced below as a source term $\mathcal{J}$ for applications with porous walls. So far no approximation has been made beyond axisymmetry and the boundary condition.

\begin{figure}[h!]
    \centering
    \includegraphics[scale=0.2]{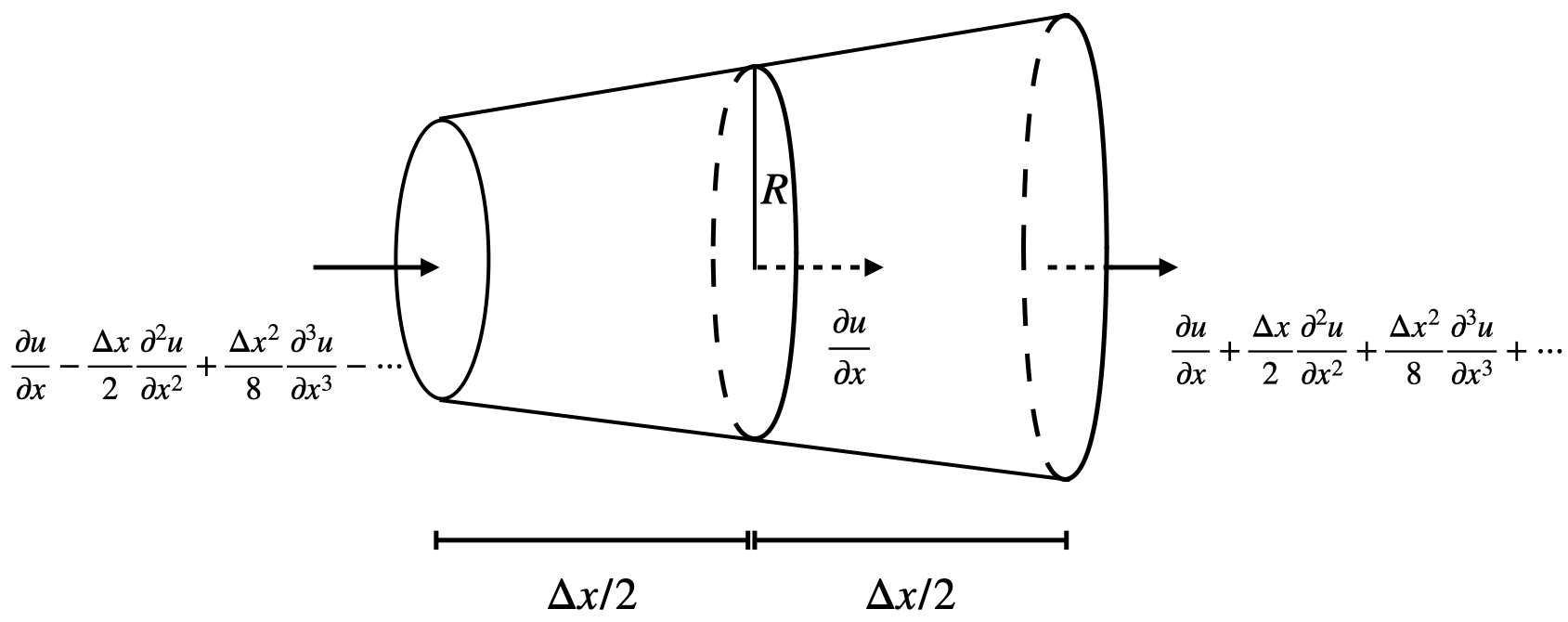}
    \caption{Expanded Jacobs flux through a truncated cone.\textbf{}}
    \label{fig:flux-expansion}
\end{figure}

To close \eqref{eq:slab-conservation} for $\bar u$, we adopt the cross-sectional uniformity ansatz $u(x, r, t) \approx \bar u(x, t)$, which is the closure shared by every reduction in the Fick-Jacobs family \cite{Jacobs1935, Zwanzig1992, RR01, KP06, DY2014}. Under the ansatz, the cap flux reduces to $\Phi(x) = -D_0 \pi R(x)^2 \partial_x \bar u(x, t)$, and \eqref{eq:slab-conservation} becomes (denoting $\bar u$ by $u$ from here on, since only the reduced 1D field appears below):
{\small
\begin{equation}
    \partial_t \int_{x_0 - \Delta x/2}^{x_0 + \Delta x/2} R(x)^2 u(x, t)\, dx \;=\; D_0 \big[ R(x_0+\Delta x/2)^2 \partial_x u(x_0+\Delta x/2, t) - R(x_0 - \Delta x/2)^2 \partial_x u(x_0 - \Delta x/2, t) \big] + \Phi_{\mathrm{lat}}/\pi. \label{eq:closed-slab}
\end{equation}
}
Equation \eqref{eq:closed-slab} is exact under the uniformity ansatz. The reduction proceeds by expanding each side in $\Delta x$ around $x_0$. The two cap fluxes on the right-hand side are the natural targets: each is the product of a cross-sectional area and an axial gradient, both smooth in $x$ and both Taylor-expanded around $x_0$. This refines Jacobs' original flux-balance derivation, which evaluated these factors at the subdomain boundaries via a first-order approximation on an infinitesimal subdomain; here, the subdomain has finite width $\Delta x$ and higher-order terms are retained. The flux configuration is shown in Fig.~\ref{fig:flux-expansion}.

Expanding the cap flux through the left face on the truncated cone (so $R(x_0 \pm \Delta x/2) = R \pm (\Delta x/2) R_x$ exactly):
\begin{displaymath}
-D_0 \left( \pi \left( R - \frac{\Delta x}{2}R_x \right)^2 \right) \left( \frac{\partial u}{\partial x} - \frac{\Delta x}{2} \frac{\partial^2 u}{\partial x^2} + \frac{\Delta x ^2}{8} \frac{\partial^3 u}{\partial x^3} -  \frac{\Delta x ^3}{48} \frac{\partial^4 u}{\partial x^4} + \frac{\Delta x^4}{384} \frac{\partial^5 u}{\partial x^5} - \cdots \right),
\end{displaymath}
and the cap flux through the right face:
\begin{displaymath}
-D_0 \left( \pi \left( R + \frac{\Delta x}{2}R_x\right)^2 \right) \left( \frac{\partial u}{\partial x} + \frac{\Delta x}{2} \frac{\partial^2 u}{\partial x^2} + \frac{\Delta x ^2}{8} \frac{\partial^3 u}{\partial x^3} +  \frac{\Delta x ^3}{48} \frac{\partial^4 u}{\partial x^4} + \frac{\Delta x^4}{384} \frac{\partial^5 u}{\partial x^5} + \cdots \right).
\end{displaymath}
The difference of cap fluxes (the right-hand side of \eqref{eq:closed-slab} apart from $\Phi_{\mathrm{lat}}/\pi$), after multiplying out and collecting by derivative order of $u$, is:
\begin{eqnarray}
    D_0 \pi \Delta x R^2 \left[ \frac{\partial^2 u}{\partial x^2} + 2 \frac{R_x}{R} \frac{\partial u}{\partial x} + \frac{\Delta x^2}{4} \frac{R_x^2}{R^2} \frac{\partial^2 u}{\partial x^2} + \frac{\Delta x^2}{4} \frac{R_x}{R} \frac{\partial^3 u}{\partial x^3} +  \frac{\Delta x^2}{24} \frac{\partial^4 u}{\partial x^4} + \frac{\Delta x^4}{96} \frac{R_x^2}{R^2} \frac{\partial^4 u}{\partial x^4} + \frac{\Delta x^4}{192} \frac{R_x}{R} \frac{\partial^5 u}{\partial x^5} \right]
\end{eqnarray}
The leading group $\partial^2 u/\partial x^2 + 2(R_x/R)\partial u/\partial x$ is $R^{-2}\partial_x(R^2 \partial_x u)$ on the cone --- the classical Fick-Jacobs flux divergence. Gathering the $\Delta x^2$ corrections in divergence form and dropping the $\Delta x^4$ tail,
\begin{equation}
       D_0 \pi \Delta x \left\{ \frac{\partial}{\partial x} \left[ \left( R^2 + \frac{\Delta x^2 R_x^2}{4} \right) \frac{\partial u}{\partial x} + \frac{\Delta x^2 R^2}{8} \frac{\partial^3 u}{\partial x^3} \right] - R^2 \frac{\Delta x^2}{12} \frac{\partial^4 u}{\partial x^4} \right\}
       \label{eq:exp-rhs}
\end{equation}
which is the $\mathcal{O}(\Delta x^2)$ right-hand side, separated into the divergence-form geometric corrections and an outer fourth-derivative residual.

For the left-hand side of \eqref{eq:closed-slab}, we use a lower order approximation as Jacobs did, approximating the subdomain integral by $\bar V \partial_t u(x_0, t)$, where $\bar V$ is the truncated-cone volume --- treating $u$ as approximately uniform across the subdomain, an additional simplification beyond the cross-sectional uniformity ansatz. Explicitly,
\begin{displaymath}
\frac{\Delta x}{3}\pi \left[ \left( R - \frac{\Delta x}{2} R_x\right)^2 + \left( R + \frac{\Delta x}{2} R_x\right)^2 + \left( R - \frac{\Delta x}{2} R_x\right) \left( R + \frac{\Delta x}{2} R_x\right) \right] \frac{\partial u}{\partial t},
\end{displaymath}
which simplifies to
\begin{equation}
    \Delta x \pi \left[ R^2 + \frac{\Delta x^2}{12} R_x^2 \right] \frac{\partial u}{\partial t}.
    \label{eq:lhs-coneVol}
\end{equation}
The bracket is the truncated-cone volume divided by $\pi \Delta x$, with the cylinder term $R^2$ and the cone-over-cylinder correction $\Delta x^2 R_x^2/12$.

No expansion is needed to restore the lateral flux $\Phi_{\mathrm{lat}}$ set aside after \eqref{eq:slab-conservation}, as it is treated over the surface of the entire element. We take the locally averaged flux $\mathcal{J}$ and multiply by the lateral surface area $\mathcal{S}$, so $\Phi_{\mathrm{lat}} = \mathcal{S}\mathcal{J}$. Setting the terms equal and dividing by $\pi \Delta x R^2$ gives the expanded divergence form of the Fick-Jacobs equation,
\begin{eqnarray}
     \left[ 1 + \frac{\Delta x^2}{12} \frac{R_x^2}{R^2} \right] \frac{\partial u}{\partial t} &=& \frac{D_0}{R^2}  \frac{\partial}{\partial x} \left[ \left( R^2 + \frac{\Delta x^2 R_x^2}{4} \right) \frac{\partial u}{\partial x} + \frac{\Delta x^2 R^2}{8} \frac{\partial^3 u}{\partial x^3} \right] - D_0 \frac{\Delta x^2}{12} \frac{\partial^4 u}{\partial x^4} \\
    &\quad& +  \frac{\mathcal{S}}{\pi\Delta x R^2} \mathcal{J}. \nonumber
\end{eqnarray}
We now have a second-order term on the right-hand-side outside of our divergence operator. This is notably geometry independent, in contrast to Eq.~\eqref{eq:fick-jacobs-exp} and of the same order as the geometric terms absorbed into the divergence operator. Thus, we do not lose model fidelity if we consider it an error term. Comparing to the $\Delta x^4$ terms dropped in Eq.~\eqref{eq:exp-rhs}, we also see that we must be in a regime such that $R_x^2/R^2 \sim 1/\Delta x^2$ for the model error to become strongly geometry dependent, a significantly less restrictive regime than $|R_x| > 1$. Then we drop the outer $\Delta x^2$ term to preserve a geometry independent model error, as well as maintain the energy stability of the model which will be demonstrated, leaving
\begin{eqnarray}
     \left[ 1 + \frac{\Delta x^2}{12} \frac{R_x^2}{R^2} \right] \frac{\partial u}{\partial t} &=& \frac{D_0}{R^2}  \frac{\partial}{\partial x} \left[ \left( R^2 + \frac{\Delta x^2 R_x^2}{4} \right) \frac{\partial u}{\partial x} + \frac{\Delta x^2 R^2}{8} \frac{\partial^3 u}{\partial x^3} \right] +  \frac{\mathcal{S}}{\pi\Delta x R^2} \mathcal{J}.
\end{eqnarray}

Notice that as $R_x, \Delta x \to 0$, this model reduces to a simple diffusion equation. Since the Fick-Jacobs model is most aptly used in regions of strong radial gradient, we compensate for that by defining a blended model. We emphasize at the outset that the blending weight $w(x)$ introduced below is a \emph{modeling construct}, not a quantity derived from the reduction above. Its role is to interpolate smoothly between pure axial diffusion in regions of negligible radial gradient and the expanded Fick-Jacobs equation in regions of strong radial gradient, two regimes for which the underlying derivation is justified. While principled construction of $w$ on specific geometries from the 3D physics is the subject of separate ongoing work, the functional form below is a specific modeling choice with the required limiting behavior. Define the weighting function
\begin{equation}
    w(x) = \frac{R_x^2}{R_x^2 + w_\delta^2}.
\end{equation}
where $w_\delta$ is a geometry-dependent parameter that defines a radial derivative threshold of when to apply the expansion model and when to apply pure diffusion. The function is chosen to be a smooth, monotone function that is bounded in the range $[0, 1]$, with the specific form chosen for computational simplicity. Then we write our complete equation
\begin{eqnarray}
     \left[ 1 + \frac{\Delta x^2}{12} \frac{R_x^2}{R^2} \right] \frac{\partial u}{\partial t} &=& w(x)\frac{D_0}{R^2}  \frac{\partial}{\partial x} \left[ \left( R^2 + \frac{\Delta x^2 R_x^2}{4} \right) \frac{\partial u}{\partial x} + \frac{\Delta x^2 R^2}{8} \frac{\partial^3 u}{\partial x^3} \right]  \nonumber \\
     &\quad& + [1-w(x)] D_0 \frac{\partial^2 u}{\partial x^2} + \frac{\mathcal{S}}{\pi\Delta x R^2} \mathcal{J}.
\end{eqnarray}
As $w(x) \to 1$ in regions of large radial derivative, the equation becomes the fully expanded Fick-Jacobs equation, and as $w(x) \to 0$ in regions of small radial derivative, the equation naturally reduces to the diffusion equation. Defining
\begin{eqnarray}
    \mathcal{A}(x) &=& 1 + \frac{\Delta x^2}{12} \frac{R_x^2}{R^2} \\
    \mathcal{F}(u) &=& \left( R^2 + \frac{\Delta x^2 R_x^2}{4} \right) \frac{\partial u}{\partial x} + \frac{\Delta x^2 R^2}{8} \frac{\partial^3 u}{\partial x^3}
\end{eqnarray}
our complete equation can be written simply as
\begin{equation}
    \mathcal{A}(x)\frac{\partial u}{\partial t} = w(x) \frac{D_0}{R^2} \frac{\partial}{\partial x} \mathcal{F}(u) + [1 - w(x)] D_0 \frac{\partial^2 u}{\partial x^2} + \frac{\mathcal{S}}{\pi\Delta x R^2} \mathcal{J} \label{eq:dFJ}
\end{equation}

\begin{remark} 
    We set $w(x) = 0$ at branch points because the assumptions of the Taylor expansion are invalid at these locations, and we may treat branch nodes as diffusion-dominated mixing regions.
\end{remark}

\begin{remark} 
    This equation naturally extends to annular regions with outer radius $R$ and inner radius $r$ simply by replacing $R^2$ with $R^2 - r^2$ and $R_x^2$ with $R_x^2 - r_x^2$.
\end{remark}

\begin{remark}
    The surface area $\mathcal{S}$ contains a factor of $\Delta x$, which will cancel with the $\Delta x$ in the denominator. Therefore, the limit $\Delta x \to 0$ is still valid, and the coefficient on surface flux terms reduces to $2/R$. 
\end{remark}



\section{Numerical Methods}\label{sec:discretization}

Define
\begin{eqnarray}
    \Delta x_{i+1/2} &:=& x_{i+1} - x_{i} \\
    \Delta x_i &:= &(x_{i+1} - x_{i-1})/2
\end{eqnarray}
Then we can define spatial derivatives
\begin{eqnarray}
    (u_{i+1/2})_x &=& \frac{u_{i+1} - u_{i}}{\Delta x_{i+1/2}} \\
    (u_i)_{xx} &=& \frac{(u_{i+1/2})_x - (u_{i-1/2})_x}{\Delta x_i} \label{eq:diff} \\
    (u_{i+1/2})_{xxx} &=& \frac{(u_{i+1})_{xx} - (u_{i})_{xx}}{\Delta x_{i+1/2}} \\
    (R_i)_x &=& \frac{1}{2} \left[ \left( \frac{R_{i+1} - R_{i}}{\Delta x_{i+1/2}} \right) + \left( \frac{R_{i} - R_{i-1}}{\Delta x_{i-1/2}} \right) \right]
\end{eqnarray}
and the semi-discrete form of Eq.~\eqref{eq:dFJ} may be written as a finite volume method,
\begin{equation}
    A_i \frac{du_i}{dt} = w_i \frac{D_0}{R_i^2} \frac{F_{i+1/2} - F_{i-1/2}}{\Delta x_i} + (1 - w_i) D_0 (u_i)_{xx} + \frac{S_i}{\pi \Delta x_i R_i^2} J_i,
\end{equation}
where $A_i$, $F_{i\pm1/2}$, $S_i$, and $J_i$ and the discrete local values of $\mathcal{A}$, $\mathcal{F}$, $\mathcal{S}$, and $\mathcal{J}$ respectively. Any appropriate time stepping method may be applied to this discretization. Specifically, we use an SBDF2 time stepping scheme \cite{BRNQ2023}, where we treat all spatial derivatives implicitly and all surface flux terms explicitly. We enforce flux balancing on the boundaries,
\begin{equation}
    F_{-1/2} = -F_{n+1/2}
\end{equation}


\subsection{Branch Nodes}

\begin{figure}
    \centering
    \includegraphics[scale=0.2]{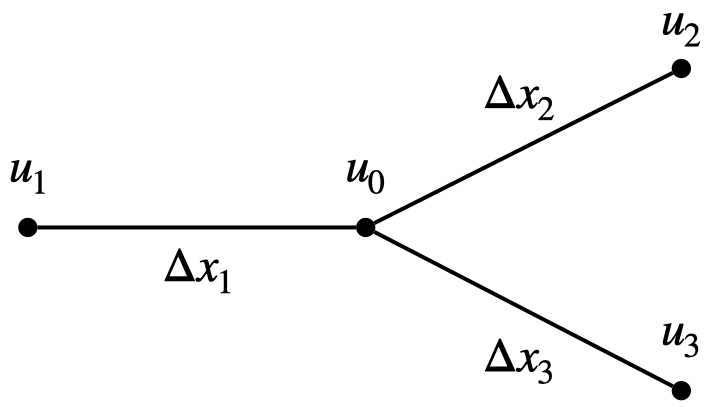}
    \caption{Discretization of a branch node.}
    \label{fig:branch}
\end{figure}

The Fick-Jacobs framework assumes slow radial variation relative to axial variation, an assumption that physically breaks down at branch nodes. In these regions, the characteristic time for radial diffusion is much shorter than the characteristic time for axial transport, resulting in nearly instantaneous radial equilibration. Therefore, we may treat branch nodes as diffusion-dominated mixing regions. In Fig.~\ref{fig:branch}, we provide a network branch with edge lengths $\Delta x_i$ and local solution values $u_i$. Then following \cite{BRNQ2023}, we define the diffusion operator,
\begin{equation}
    (u_0)_{xx} = \sum_{i > 0} \frac{u_i}{\Delta x_i} - u_0 \sum_{i > 0} \frac{1}{\Delta x_i}.
\end{equation}
The diffusion operator may be applied everywhere in this form. On a straight cable (a single node with two neighbors), this exactly reduces to Eq.~\eqref{eq:diff}.

\section{Stability and Convergence}\label{sec:stability-convergence}

In this section, we establish the stability of the discrete model, and then use this result to demonstrate convergence. We apply the following assumptions, the first being a strict definition of the truncated cone geometry:
\begin{itemize}
    \item[(A1)] \textbf{Geometry.} The radius $R \in C^4([0,L])$ is uniformly bounded below by $R_{\min} > 0$ such that $R(x) \geq R_{\min}$ for all $x \in [0, L]$. Equivalently, the cross-sectional area $A(x) = \pi R(x)^2$ is bounded away from zero, $A(x) \geq \pi R_{\min}^2 > 0$.
    \item[(A2)] \textbf{Regularity of the solution.} The solution $u \in C^{4,3}([0,L] \times [0,T])$ with $\|u\|_{C^{4,3}} \leq C_u$ for a finite constant $C_u$, where $\|u\|_{C^{4,3}} := \sum_{|\alpha| \leq 4} \|\partial_x^\alpha u\|_{L^\infty} + \sum_{\beta \leq 3} \|\partial_t^\beta u\|_{L^\infty}$.
\end{itemize}

\begin{lemma}[Discrete energy stability for the expansion model]\label{stab-lemma}
    Consider the expansion-model flux $\mathcal{F}(u) = a(x) \partial_x u + b(x) \partial_x^3 u$ on $[0, L]$ with
    \begin{equation}
        a(x) \;=\; R(x)^2 + \frac{\Delta x^2}{4} R_x(x)^2, \qquad b(x) \;=\; \frac{\Delta x^2}{8} R(x)^2, \label{eq:expansion-ab}
    \end{equation}
    and the finite-volume semi-discretization of $\partial_t u = (D_0/R^2)\, \partial_x \mathcal{F}(u)$ on a uniform grid of spacing $\Delta x$, with no-flux boundary conditions $F_{-1/2} = F_{n+1/2} = 0$ and numerical conditions $(u_0)_{xx} = (u_n)_{xx} = 0$. Suppose
    \begin{equation}
        \Delta x \;\leq\; \Delta x_*, \qquad \Delta x_* \;:=\; \frac{R_{\min}^2}{\|R\|_\infty^2 + \|R\|_\infty \|R_x\|_\infty}. \label{eq:stab-threshold}
    \end{equation}
    Given geometry (A1), the discretization is discretely $L^2$ energy stable in the area-weighted norm $\|u\|_R^2 := \sum_i u_i^2 R_i^2 \Delta x$: $dE/dt \leq 0$ for the homogeneous problem, where $E := \tfrac{1}{2} \|u\|_R^2$.
\end{lemma}

\begin{proof}
   The proof has three parts: setup with the area-weighted energy; summation by parts and decomposition of the third-derivative flux contribution; and a combined estimate that closes via the threshold \eqref{eq:stab-threshold}.

    \emph{Setup.}
    The full equation \eqref{eq:dFJ} with $w(x) \equiv 1$ reads $\mathcal{A}(x)\, \partial_t u = (D_0/R^2)\, \partial_x \mathcal{F}(u)$, where $\mathcal{A}(x) = 1 + (\Delta x^2/12)(R_x/R)^2$ is bounded with $1 \le \mathcal{A}(x) \le 1 + (\Delta x^2/12)(\|R_x\|_\infty/R_{\min})^2$, in particular uniformly bounded above and below. We absorb $\mathcal{A}$ into the time variable by a positive rescaling that does not affect the sign of $dE/dt$, setting $\mathcal{A} \equiv 1$ for the stability analysis without loss of generality. The semi-discrete scheme is
    \begin{equation}
        \dot{u}_i = \frac{D_0}{R_i^2}\, \frac{F_{i+1/2} - F_{i-1/2}}{\Delta x},
    \end{equation}
    with $F_{-1/2} = F_{n+1/2} = 0$. Using the area-weighted inner product $\langle u, v\rangle_R := \sum_i u_i v_i R_i^2 \Delta x$ and the corresponding energy $E := \tfrac{1}{2}\langle u, u\rangle_R$, the $R^2$ weight cancels the $D_0/R_i^2$ factor on the right-hand side:
    \begin{equation}
        \frac{dE}{dt} \;=\; \sum_i u_i\, R_i^2 \, \dot u_i\, \Delta x \;=\; D_0 \sum_i u_i (F_{i+1/2} - F_{i-1/2}).
    \end{equation}
    Using $F_{-1/2} = F_{n+1/2} = 0$ to discard boundary terms, summation by parts gives
    \begin{equation}
        \frac{dE}{dt} \;=\; D_0 \sum_i (u_i - u_{i+1})\, F_{i+1/2} \;=\; -D_0 \sum_i \Delta x\, (u_{i+1/2})_x\, F_{i+1/2},
    \end{equation}
    where $(u_{i+1/2})_x := (u_{i+1} - u_i)/\Delta x$. Substituting $F_{i+1/2} = a_{i+1/2}(u_{i+1/2})_x + b_{i+1/2}(u_{i+1/2})_{xxx}$,
    \begin{equation}
        \frac{dE}{dt} \;=\; \underbrace{-D_0 \Delta x \sum_i a_{i+1/2}\, (u_{i+1/2})_x^2}_{T_1} \;+\; \underbrace{-D_0 \Delta x \sum_i b_{i+1/2}\, (u_{i+1/2})_x\, (u_{i+1/2})_{xxx}}_{T_2}.
    \end{equation}
    The first sum is non-positive, $T_1 \leq 0$. The second sum has indefinite sign in general.

    \emph{Decomposition of $T_2$.}
    Set $D^+ u_i := (u_{i+1} - u_i)/\Delta x$ and $(u_i)_{xx} := (u_{i+1} - 2u_i + u_{i-1})/\Delta x^2$. By definition of the third-derivative finite difference,
    \begin{equation}
        (u_{i+1/2})_{xxx} \;=\; \frac{(u_{i+1})_{xx} - (u_i)_{xx}}{\Delta x}.
    \end{equation}
    Substituting and summing by parts on $T_2$ using the numerical boundary conditions $(u_0)_{xx} = (u_n)_{xx} = 0$,
    \begin{equation}
        T_2 \;=\; -\frac{D_0}{\Delta x} \sum_i b_{i+1/2}\, D^+ u_i\, \big[(u_{i+1})_{xx} - (u_i)_{xx}\big]\,\Delta x \;=\; D_0 \Delta x \sum_i \big[b_{i+1/2} D^+ u_i - b_{i-1/2} D^+ u_{i-1}\big](u_i)_{xx}.
    \end{equation}
    Splitting the bracket as $b_{i+1/2}(D^+ u_i - D^+ u_{i-1}) + (b_{i+1/2} - b_{i-1/2}) D^+ u_{i-1}$ and using $D^+ u_i - D^+ u_{i-1} = \Delta x\, (u_i)_{xx}$ on a uniform grid,
    \begin{equation}
        T_2 \;=\; \underbrace{D_0 \Delta x^2 \sum_i b_{i+1/2}\, \big((u_i)_{xx}\big)^2}_{T_2^{(a)} \,\geq\, 0} \;+\; \underbrace{D_0 \Delta x \sum_i (b_{i+1/2} - b_{i-1/2})\, D^+ u_{i-1}\, (u_i)_{xx}}_{T_2^{(b)}}. \label{eq:t2-decomp}
    \end{equation}

    \emph{Combined estimate.}
    For $T_2^{(a)}$, substitute the expansion-model $b_{i+1/2} = \Delta x^2 R_{i+1/2}^2/8$ and the inverse inequality $((u_i)_{xx})^2 \leq (2/\Delta x^2)\big((D^+ u_i)^2 + (D^+ u_{i-1})^2\big)$ (from $(u_i)_{xx} = (D^+ u_i - D^+ u_{i-1})/\Delta x$ and $(p-q)^2 \le 2(p^2 + q^2)$):
    \begin{equation}
        T_2^{(a)} \;=\; \frac{D_0 \Delta x^4}{8} \sum_i R_{i+1/2}^2 \big((u_i)_{xx}\big)^2 \;\leq\; \frac{D_0 \Delta x^2}{4} \sum_i R_{i+1/2}^2 \big[(D^+ u_i)^2 + (D^+ u_{i-1})^2\big].
    \end{equation}
    The numerical conditions $(u_0)_{xx} = (u_n)_{xx} = 0$ restrict the cell sum to interior cells $i \in \{1, \ldots, n-1\}$, on which both $D^+ u_i$ (at face $i + 1/2$) and $D^+ u_{i-1}$ (at face $i - 1/2$) are well-defined interior-face forward differences. Bounding $R_{i+1/2}^2 \leq \|R\|_\infty^2$ pointwise and applying the change of summation index $j = i-1$ to the second sum,
    \begin{equation}
        \sum_{i=1}^{n-1} R_{i+1/2}^2 \big[(D^+ u_i)^2 + (D^+ u_{i-1})^2\big] \;\leq\; \|R\|_\infty^2 \bigg[ \sum_{i=1}^{n-1} (D^+ u_i)^2 + \sum_{j=0}^{n-2} (D^+ u_j)^2 \bigg].
    \end{equation}
    Both bracketed sums consist of non-negative terms and are subsets of the interior-face sum $S := \sum_{k=0}^{n-1}(D^+ u_k)^2$, so each is bounded by $S$. Therefore
    \begin{equation}
        T_2^{(a)} \;\leq\; \frac{D_0 \Delta x^2}{2}\, \|R\|_\infty^2\, S. \label{eq:t2a-bound-clean}
    \end{equation}
    The lower bound $|T_1| = D_0 \Delta x \sum_{k=0}^{n-1} a_{k+1/2}(D^+ u_k)^2 \geq D_0 \Delta x\, R_{\min}^2\, S$ (using $a \geq R^2 \geq R_{\min}^2$) gives $S \leq |T_1|/(D_0 \Delta x\, R_{\min}^2)$, and
    \begin{equation}
        T_2^{(a)} \;\leq\; \frac{\Delta x\, \|R\|_\infty^2}{2 R_{\min}^2}\, |T_1|.
    \end{equation}

    For $T_2^{(b)}$, $b_{i+1/2} - b_{i-1/2} = (\Delta x^2/8)(R_{i+1/2}^2 - R_{i-1/2}^2)$, so on a uniform grid the mean value theorem (applied to $R^2$ over $[x_{i-1/2}, x_{i+1/2}]$) gives $|R_{i+1/2}^2 - R_{i-1/2}^2| \leq 2 \|R R_x\|_\infty \Delta x \leq 2\|R\|_\infty \|R_x\|_\infty \Delta x$, hence
    \begin{equation}
        |b_{i+1/2} - b_{i-1/2}| \;\leq\; \frac{\Delta x^3}{4}\, \|R\|_\infty\, \|R_x\|_\infty.
    \end{equation}
    Then
    \begin{equation}
        |T_2^{(b)}| \;\leq\; \frac{D_0 \Delta x^4}{4}\, \|R\|_\infty\, \|R_x\|_\infty \sum_i |D^+ u_{i-1}|\, |(u_i)_{xx}|.
    \end{equation}
    Apply Cauchy--Schwarz to $|T_2^{(b)}|$:
    \begin{equation}
        \sum_i |D^+ u_{i-1}|\, |(u_i)_{xx}| \;\leq\; \bigg(\sum_i (D^+ u_{i-1})^2\bigg)^{1/2} \bigg(\sum_i ((u_i)_{xx})^2\bigg)^{1/2}.
    \end{equation}
    The first factor is bounded by $\sqrt{S}$ (by the change of summation index $j = i-1$, the sum equals $\sum_{j=0}^{n-2}(D^+ u_j)^2 \leq S$). The second factor is bounded via the inverse inequality used for $T_2^{(a)}$,
    \begin{equation}
        \sum_i ((u_i)_{xx})^2 \;\leq\; \frac{2}{\Delta x^2} \bigg[\sum_i (D^+ u_i)^2 + \sum_i (D^+ u_{i-1})^2\bigg] \;\leq\; \frac{4}{\Delta x^2}\, S,
    \end{equation}
    whose square root is $(2/\Delta x)\sqrt{S}$. Combining,
    \begin{equation}
        |T_2^{(b)}| \;\leq\; \frac{D_0 \Delta x^4}{4}\, \|R\|_\infty\, \|R_x\|_\infty \cdot \sqrt{S} \cdot \frac{2\sqrt{S}}{\Delta x} \;=\; \frac{D_0 \Delta x^3}{2}\, \|R\|_\infty\, \|R_x\|_\infty\, S \;\leq\; \frac{\Delta x^2}{2}\, \frac{\|R\|_\infty \|R_x\|_\infty}{R_{\min}^2}\, |T_1|.
    \end{equation}

    \emph{Conclusion.} Combining the two bounds,
    \begin{equation}
        T_2 \;\leq\; T_2^{(a)} + |T_2^{(b)}| \;\leq\; \frac{\Delta x}{2 R_{\min}^2} \left[ \|R\|_\infty^2 + \Delta x\, \|R\|_\infty \|R_x\|_\infty \right] |T_1|.
    \end{equation}
    For the bracketed coefficient on $|T_1|$ to be at most $1$, it suffices that \newline $\Delta x \|R\|_\infty^2/(2 R_{\min}^2) + \Delta x^2 \|R\|_\infty \|R_x\|_\infty/(2 R_{\min}^2) \leq 1$. A sufficient (and slightly stronger) condition is $\Delta x (\|R\|_\infty^2 + \|R\|_\infty \|R_x\|_\infty)/R_{\min}^2 \leq 1$, which is exactly the threshold \eqref{eq:stab-threshold}. Under this condition, $T_1 + T_2 \leq 0$, giving $dE/dt \leq 0$.
\end{proof}

\begin{remark}[Boundary conditions]
    The numerical boundary conditions $(u_0)_{xx} = (u_n)_{xx} = 0$ are imposed only on the $\Delta x^2$ expansion terms of $\mathcal{F}(u)$, which vanish under mesh refinement and therefore do not affect the physical solution.
\end{remark}

\begin{remark}[Verification on the test geometry]
    The manufactured-solution test geometry in Sec.~\ref{sec:numerical-results} has $R(x) = 1 + 2\sin(x)$ on the domain bounded by the first two positive roots of $2\cos(2x) - \sin(x) = 0$, giving $\|R_x\|_\infty \approx 1.61$, $R_{\min} \approx 2.19$, and $\|R\|_\infty \approx 2.98$. The threshold \eqref{eq:stab-threshold} requires $\Delta x \leq R_{\min}^2/(\|R\|_\infty^2 + \|R\|_\infty\|R_x\|_\infty) \approx 4.80/(8.88 + 4.80) \approx 0.35$, comfortably satisfied by the grid spacings used in the numerical experiments. The threshold is a small-$\Delta x$ condition with constants depending on the geometry; it tightens as $R_{\min}$ shrinks (narrow constrictions) and loosens for tubes of uniform radius. The convergence regime $\Delta x \to 0$ relevant for Theorem~\ref{thm:convergence} automatically satisfies the threshold for any fixed smooth geometry.
\end{remark}


Now that we have established stability, we only need to demonstrate consistency to establish convergence. We show this below.

\begin{theorem}[Convergence with explicit rate]\label{thm:convergence}
Suppose assumptions (A1)--(A2) hold and that $w(x) \equiv 1$ on $[0, L]$. Given initial data $u_0 \in C^4([0, L])$ and no-flux boundary conditions $\partial_x u(0, t) = \partial_x u(L, t) = 0$ for $t \in [0, T]$, let $u(x, t)$ denote the analytical solution of the Fick-Jacobs equation \eqref{eq:fick-jacobs}, and let $u^h_i(t)$ denote the semi-discrete solution of the finite-volume discretization \eqref{eq:dFJ} on a uniform grid $\{x_i\}_{i=0}^n$ of spacing $\Delta x$ with the same initial data $u^h_i(0) = u_0(x_i)$ and no-flux discrete boundary conditions $F_{-1/2} = F_{n+1/2} = 0$. Define the discrete $L^2$ error $\|e(t)\|_h := \big( \sum_i (u(x_i, t) - u^h_i(t))^2 \, \Delta x \big)^{1/2}$. Then there exists a constant $C = C(C_u, R_{\min}, \|R\|_{C^3}, T, D_0)$, independent of $\Delta x$, such that for all $t \in [0, T]$,
    \begin{equation}
        \|e(t)\|_h \leq C \, \Delta x^2. \label{eq:convergence-rate}
    \end{equation}
    The same conclusion holds for the fully discrete SBDF2 scheme with combined error bound $\mathcal{O}(\Delta x^2 + \Delta t^2)$, under the standard SBDF2 stability condition on $\Delta t$.
\end{theorem}

\begin{proof}
    The argument has three parts: identification of the limiting model under $\Delta x \to 0$, consistency of the discretization with the limiting model at rate $\mathcal{O}(\Delta x^2)$, and propagation via the discrete stability estimate of Lemma~\ref{stab-lemma} to the global error bound.

    \emph{Step 1: identification of the limiting model.}
    Given a smooth domain $\Omega$ with boundary $\Gamma$. Define operators
    \begin{eqnarray}
        \mathcal{L}_L u &=& w(x) \frac{D_0}{R^2} \frac{\partial}{\partial x} \left[ \left( R^2 + \frac{\Delta x^2 R_x^2}{4} \right) \frac{\partial u}{\partial x} \right] ,\\
        \mathcal{L}_H u &=& w(x) \frac{D_0}{R^2} \frac{\partial}{\partial x} \left[\frac{\Delta x^2 R^2}{8} \frac{\partial^3 u}{\partial x^3} \right], \\
        \mathcal{L}_D u &=& [1 - w(x)] D_0 \frac{\partial^2 u}{\partial x^2}.
    \end{eqnarray}
    Assuming $w(x)=1$, this leaves
    \begin{eqnarray}
        \mathcal{L}_L u &=& \frac{D_0}{R^2} \frac{\partial}{\partial x} \left[ \left( R^2 + \frac{\Delta x^2 R_x^2}{4} \right) \frac{\partial u}{\partial x} \right], \\
        \mathcal{L}_H u &=& \frac{D_0}{R^2} \frac{\partial}{\partial x} \left[\frac{\Delta x^2 R^2}{8} \frac{\partial^3 u}{\partial x^3} \right].
    \end{eqnarray}
    Apply no-flux boundary conditions to each operator, and the boundary condition
    \begin{equation}
        \left. \frac{\partial^2 u}{\partial x^2} \right|_{x \in \Gamma} = 0
    \end{equation}
    to the operator $\mathcal{L}_H$. The coefficients in $\mathcal{L}_L$ and $\mathcal{L}_H$ depend explicitly on $\Delta x$. In the limit $\Delta x \to 0$, the operator $\mathcal{L}_H \to 0$ pointwise (the coefficient of $\partial^3_x u$ vanishes), and $\mathcal{L}_L$ reduces to the standard Fick-Jacobs operator $(D_0/R^2)\partial_x(R^2 \partial_x u)$. The volumetric coefficient $\mathcal{A}(x) = 1 + (\Delta x^2/12)(R_x^2/R^2)$ similarly reduces to unity. Thus the continuous limit of the modified PDE is the Fick-Jacobs equation \eqref{eq:fick-jacobs}.

    \emph{Step 2: consistency.}
    Define discrete derivative approximations,
    \begin{eqnarray}
        (u_{i+1/2})_x &=& \frac{u_{i+1} - u_{i}}{\Delta x} ,\\
        (u_i)_{xx} &=& \frac{u_{i+1} - 2u_i + u_{i-1}}{\Delta x^2}.
    \end{eqnarray}
    Let $L_L$ and $L_H$ denote the corresponding discrete operators acting at cell $i$. We derive the truncation bound for $L_L$ in detail, with $L_H$ analogous.

    The discrete operator at cell $i$ is $L_L u_i = (D_0/R_i^2)\big[a_{i+1/2}(D^+ u)_i - a_{i-1/2}(D^+ u)_{i-1}\big]/\Delta x$. Taylor-expanding $u$ around $x_{i+1/2}$ to fourth order under (A2),
    \begin{equation}
        (D^+ u)_i \;=\; \frac{u(x_{i+1/2} + \Delta x/2) - u(x_{i+1/2} - \Delta x/2)}{\Delta x} \;=\; \partial_x u(x_{i+1/2}) + \frac{\Delta x^2}{24}\, \partial_x^3 u(x_{i+1/2}) + O(\Delta x^4 \, C_u),
    \end{equation}
    so $a_{i+1/2}(D^+ u)_i = g(x_{i+1/2}) + (\Delta x^2/24)\, a(x_{i+1/2})\, \partial_x^3 u(x_{i+1/2}) + O(\Delta x^4 \, C_u\, \|a\|_\infty)$, where $g(x) := a(x)\partial_x u(x)$. The same expansion at $x_{i-1/2}$ gives the analogous expression. Subtracting and dividing by $\Delta x$,
    {\small
    \begin{equation}
        \frac{a_{i+1/2}(D^+ u)_i - a_{i-1/2}(D^+ u)_{i-1}}{\Delta x} \;=\; \frac{g(x_{i+1/2}) - g(x_{i-1/2})}{\Delta x} + \frac{\Delta x^2}{24} \, \frac{a\,\partial_x^3 u\big|_{x_{i+1/2}} - a\,\partial_x^3 u\big|_{x_{i-1/2}}}{\Delta x} + O(\Delta x^3 \, C_u\, \|a\|_\infty).
    \end{equation}
    }
    Each of the two quotients on the right is a central difference. Under (A2) and $g \in C^3$ (since $a \in C^3$ ($R \in C^4$ via (A1)) and $\partial_x u \in C^3$ ($u \in C^4$ via (A2)), the first central difference equals $\partial_x g(x_i) + O(\Delta x^2 \|g\|_{C^3}) = \partial_x[a\partial_x u]\big|_{x_i} + O(\Delta x^2 C_u \|R\|_{C^4}^2)$. The second central difference is $O(1)$ in $\Delta x$, so the full contribution is $O(\Delta x^2 C_u \|a\|_\infty) = O(\Delta x^2 C_u \|R\|^2_\infty)$. Multiplying by $D_0/R_i^2 \le D_0/R_{\min}^2$,
    \begin{equation}
        L_L u(x_i) \;=\; \mathcal{L}_L u(x_i) + \tau^L_i, \qquad |\tau^L_i| \;\leq\; C_L\, \Delta x^2\, C_u\, \|R\|_{C^4}/R_{\min}^2, \label{eq:LL-consistency}
    \end{equation}
    where $C_L$ is an absolute constant collecting numerical factors.

    For $L_H$, the same analysis applies to the flux $b(x)\partial_x^3 u$. The discrete third derivative $(u_{i+1/2})_{xxx} = [(u_{i+1})_{xx} - (u_i)_{xx}]/\Delta x$ has truncation error $O(\Delta x \, C_u)$ under (A2), since $u \in C^4$ provides only one half-order of accuracy beyond the third-derivative approximation. However, $b(x) = \Delta x^2 R^2/8$ carries an explicit $\Delta x^2$ prefactor, so the truncation error of $L_H$ is $O(\Delta x^2 \cdot \Delta x \, C_u\, \|R\|^2_\infty / R_{\min}^2)$, which is $O(\Delta x^3)$ — \emph{sub-leading} relative to the $L_L$ contribution. Writing the dominant $O(\Delta x^2)$ part explicitly,
    \begin{equation}
        L_H u(x_i) \;=\; \mathcal{L}_H u(x_i) + \tau^H_i, \qquad |\tau^H_i| \;\leq\; C_H\, \Delta x^2\, C_u\, \|R\|_{C^4}/R_{\min}^2,
    \end{equation}
    which is the bound used below; the actual order is sharper but is not needed for the headline rate.

    The dominant contribution to the combined truncation comes from $L_L$, and the rate is $\mathcal{O}(\Delta x^2)$. The $C^4$ regularity in (A2) is sharp for this rate.

    \emph{Step 3: propagation to the global error.}
    Let $e_i(t) := u(x_i, t) - u^h_i(t)$ denote the pointwise error, where $u$ is the analytical FJ solution and $u^h_i$ is the discrete expansion-model solution. Subtracting the semi-discrete scheme from the analytical FJ equation evaluated at grid points and adding/subtracting $(L_L + L_H) u(x_i)$ gives
    \begin{equation}
        \frac{d e_i}{dt} = (L_L + L_H) e_i + \tau_i(t), \qquad \tau_i := \underbrace{[\mathcal{L}_{\rm FJ} - \mathcal{L}_L - \mathcal{L}_H] u(x_i)}_{\delta^{\rm mod}_i} \;-\; \tau^L_i \;-\; \tau^H_i,
    \end{equation}
    where $\mathcal{L}_{\rm FJ} u := (D_0/R^2)\partial_x(R^2 \partial_x u)$ is the standard FJ operator. The term $\delta^{\rm mod}_i$ is the \emph{modeling residual}, measuring the difference between the FJ operator and the modified operator $\mathcal{L}_L + \mathcal{L}_H$ applied to the same function $u$. Direct calculation gives $\mathcal{L}_{\rm FJ} - \mathcal{L}_L - \mathcal{L}_H = -(D_0/R^2) \partial_x[(\Delta x^2 R_x^2/4)\partial_x u + (\Delta x^2 R^2/8) \partial_x^3 u]$, so under (A1)--(A2),
    \begin{equation}
        |\delta^{\rm mod}_i| \leq C_{\rm mod} \, \Delta x^2 \, C_u \, \|R\|_{C^4} / R_{\min}^2.
    \end{equation}
    The terms $\tau^L_i$ and $\tau^H_i$ are the discretization truncation errors bounded in Step~2. Each of the three contributions is $\mathcal{O}(\Delta x^2)$ uniformly in $\Delta x$, so $|\tau_i(t)| = \mathcal{O}(\Delta x^2)$ uniformly in $i$ and $t$. Lemma~\ref{stab-lemma} establishes that the discrete operator $L_L + L_H$ acting on the error $e$ satisfies $\langle e, (L_L + L_H) e \rangle_R \leq 0$ in the area-weighted discrete inner product $\langle u, v\rangle_R := \sum_i u_i v_i R_i^2 \Delta x_i$. Under (A1), $\langle\cdot, \cdot\rangle_R$ and the unweighted $\langle u, v\rangle_h := \sum_i u_i v_i \Delta x_i$ induce equivalent norms with ratio bounded by $\|R\|_\infty^2/R_{\min}^2$. Taking the inner product of the error equation with $e$ in the weighted inner product gives
    \begin{equation}
        \frac{1}{2} \frac{d}{dt} \|e\|_R^2 \;=\; \langle e, (L_L + L_H) e\rangle_R + \langle e, \tau\rangle_R \;\leq\; \|e\|_R \cdot \|\tau\|_R,
    \end{equation}
    where $\|e\|_R^2 := \langle e, e\rangle_R$, the first term is non-positive by Lemma~\ref{stab-lemma}, and the second inequality uses Cauchy--Schwarz. Standard discrete Gronwall (see, e.g., \cite{ThomasTAM22}) applied to this differential inequality gives
    \begin{equation}
        \|e(t)\|_R \;\leq\; \|e(0)\|_R + \int_0^t \|\tau(s)\|_R\, ds \;\leq\; C T\, \Delta x^2\, (1 + C_u) \label{eq:gronwall-bound}
    \end{equation}
    assuming the initial discretization error vanishes (e.g., by $u^h_i(0) = u(x_i, 0)$). Since $|\tau_i| = \mathcal{O}(\Delta x^2)$ uniformly, $\|\tau\|_R \leq C(L, \|R\|_\infty)\Delta x^2$. By norm equivalence (under (A1)), the unweighted discrete $L^2$ bound $\|e(t)\|_h \le C \Delta x^2$ stated in \eqref{eq:convergence-rate} follows with a constant inflated by at most $\|R\|_\infty/R_{\min}$. The discrete $L^1$ error inherits the same rate by Cauchy--Schwarz on $[0,L]$.

    The fully discrete SBDF2 statement follows from standard analysis of IMEX schemes \cite{ARW1995}: under (A2), time discretization adds $\mathcal{O}(\Delta t^2)$ to the truncation error, with additive propagation under the unconditional implicit-diffusion stability.
\end{proof}


\subsection{On the singular-perturbation structure and the scope of Theorem~\ref{thm:convergence}}\label{sec:singular-perturbation}

The discretization parameter $\Delta x$ enters the model in two distinct places: as the grid spacing of the finite-volume discretization, and as a parameter in the coefficients of the modified continuous PDE \eqref{eq:dFJ} (in $\mathcal{A}(x)$ and $\mathcal{F}(u)$). In implementation, these are identified and go to zero together under mesh refinement. This is not the conventional Lax setting of a fixed PDE with refining grid, and the convergence statement of Theorem~\ref{thm:convergence} requires an argument accommodating the singular-perturbation structure.

\subsubsection{Decomposition of the error}\label{sec:error-decomp}

Let $u^{\rm FJ}(x, t)$ denote the analytical solution of the classical Fick-Jacobs equation \eqref{eq:fick-jacobs}, $u^{(\Delta x)}(x, t)$ the analytical solution of the modified PDE \eqref{eq:dFJ} at the given $\Delta x$ (with $w(x) \equiv 1$ and matching initial/boundary data), and $u^h_i(t)$ the discrete solution on the grid of spacing $\Delta x$. The error admits the triangle decomposition
\begin{equation}
    \big\| u^{\rm FJ}(x_i, t) - u^h_i(t) \big\|_h \;\leq\; \underbrace{\big\| u^{\rm FJ}(x_i, t) - u^{(\Delta x)}(x_i, t) \big\|_h}_{\text{(I) modeling error}} \;+\; \underbrace{\big\| u^{(\Delta x)}(x_i, t) - u^h_i(t) \big\|_h}_{\text{(II) numerical error}}. \label{eq:triangle-decomp}
\end{equation}

\textbf{Term (I), the modeling error,} measures the difference between the analytical FJ and modified-PDE solutions at the same $\Delta x$. The modified PDE differs from FJ by $\mathcal{O}(\Delta x^2)$ perturbations in its coefficients, so standard continuous dependence of parabolic solutions on coefficients \cite{ThomasTAM22} gives $\| u^{\rm FJ}(\cdot, t) - u^{(\Delta x)}(\cdot, t) \|_{L^2([0,L])} \leq C_1(t)\, \Delta x^2$, with $C_1$ independent of $\Delta x$.

\textbf{Term (II), the numerical error,} is the classical Lax-equivalence question, with the subtlety that the PDE depends on $\Delta x$. The consistency analysis in Step~2 of the proof of Theorem~\ref{thm:convergence} establishes that the discrete operators $L_L$, $L_H$ are consistent with their continuous counterparts at rate $\mathcal{O}(\Delta x^2)$ uniformly in $\Delta x$. The stability estimate of Lemma~\ref{stab-lemma} holds uniformly in $\Delta x$ (the coefficients $a, b$ are bounded with bounded derivatives under (A1)). Combining gives $\| u^{(\Delta x)}(x_i, t) - u^h_i(t) \|_h \leq C_2(t)\, \Delta x^2$, with $C_2$ depending on $T$, the regularity of $u^{(\Delta x)}$ (uniform in $\Delta x$; see Remark~\ref{rem:uniform-regularity}), and $D_0$, but not on $\Delta x$.

Summing the two bounds (and absorbing the $L^2 \to$ discrete-$L^2$ adjustment) gives the $\mathcal{O}(\Delta x^2)$ bound stated in Theorem~\ref{thm:convergence}.

Two further points clarify the scope. Theorem~\ref{thm:convergence} certifies that the discretization is a faithful numerical approximation of the modified PDE \eqref{eq:dFJ} at each $\Delta x$ (term II), and that the discrete solution converges to the FJ solution in the joint limit (the full I+II bound). It does \emph{not} establish modeling fidelity at finite $\Delta x$ --- that the modified PDE is a closer reduction of the 3D axisymmetric diffusion equation than FJ is at the same $\Delta x$. That is a statement about \emph{model error} rather than numerical error, tested directly in the manufactured-solution experiment of Sec.~\ref{sec:numerical-results} (radial-gradient-independence of the total error) and in the axisymmetric benchmark of Sec.~\ref{sec:bpb-benchmark} (comparison to 3D Brownian dynamics on the BPB cone after tuning of $w_\delta$).

\begin{remark}[Uniform regularity of $u^{(\Delta x)}$]\label{rem:uniform-regularity}
    The bound on term (II) requires $\|u^{(\Delta x)}\|_{C^4([0,L])}$ to remain bounded as $\Delta x \to 0$. This follows from standard parabolic regularity since \eqref{eq:dFJ} is a smooth perturbation of FJ in its coefficients and $u^{\rm FJ}$ has the required regularity by (A2).
\end{remark}

\begin{remark}[On the role of $w(x)$]
    Theorem~\ref{thm:convergence} is stated for $w(x) \equiv 1$. For $w(x) \in [0, 1]$, the model is a convex combination of the pure expansion operator and the standard diffusion operator; the diffusion contribution is independent of $\Delta x$ and lies in the conventional Lax setting, while the $w$-weighted contribution follows the argument above. Both are $L^2$-energy-stable and $\mathcal{O}(\Delta x^2)$-consistent, and convex combination preserves the rate.
\end{remark}



\section{Numerical Results}\label{sec:numerical-results}

In this section we present three numerical examples. First, an exact-solution test on a channel of varying radius demonstrates the stability, accuracy, and convergence of the method against FJ and competing corrections. Second, a direct benchmark against the three-dimensional axisymmetric diffusion equation on the truncated-cone geometry of Berezhkovskii, Pustovoit, and Bezrukov \cite{BPB2007} verifies that the expansion model --- after per-$\lambda$ tuning of $w_\delta$ --- reproduces the mean first-passage time across the full slope range, including the regime $|R_x| > 1$ where previous corrections fail. Third, a biological application compares 1D to full 3D simulation, demonstrating qualitative dynamics that the standard 1D reduction misses.


\subsection{1D Fick-Jacobs Model}

In this example, we compare our expansion model with the classical Fick-Jacobs model (FJ) and the Zwanzig (Zw), Reguera-Rubi (RR), Kalinay-Percus (KP), Kalinay (Kal), and Dorfman-Yariv (DY) corrections. Define solutions of the form
\begin{equation}
    u(x,t) = (1 + \alpha \sin(kx)) e^{-D_0k^2t}\cos(kx) + \beta,
\end{equation}
where $1 + \alpha \sin(kx)$ defines the radius of the domain and $\beta$ is a constant chosen to enforce positivity in the solution. These satisfy Eq.~\eqref{eq:fick-jacobs} with additional volumetric forcing,
\begin{equation}
    \frac{\partial u}{\partial t} = \frac{D_0}{R^2} \frac{\partial}{\partial x} \left( R^2 \frac{\partial u}{\partial x} \right) + f(x),
\end{equation}
where
\begin{equation}
    f(x) = \frac{\alpha D_0 e^{-D_0 k^2 t} k^2 \cos(kx)}{2(1 + \alpha \sin(kx))}(3\alpha - 7\alpha \cos(2kx) + 10\sin(kx)).
\end{equation}
Because this is a volumetric flux and not a surface flux, we make the scaling adjustment
\begin{equation}
    \mathcal{J} = \frac{\mathcal{V}}{\mathcal{S}} f(x)
\end{equation}
to incorporate this into our model, where $\mathcal{V}$ is the volume of the domain under consideration. Then the surface area term here cancels with the surface area term in Eq.~\eqref{eq:dFJ}, and we are left with the flux over the whole volume. We define the domain such that we maintain zero-flux ($u_x = 0$) on the boundaries. We have
\begin{equation}
    \frac{\partial u}{\partial x} = e^{D_0k^2t}k(\alpha \cos(2kx) - \sin(kx)),
\end{equation}
so we define our domain by, for simplicity, numerically solving for the first two solutions of
\begin{equation}
    \alpha \cos(2kx) - \sin(kx) = 0.
\end{equation}
For our test, we choose $\alpha = 2$ (giving $\max |R_x| = 1.61$, beyond the range of validity), $k = 1$, $D_0 = 10^{-2}$, $\beta = 2$, and set $\Delta t = 10^{-2}$ and final time $T = 500$. The initial condition fixes the unique solution for a given $\beta$ in numerical computations.

We evaluate the numerical solutions by computing the error from the exact solution. Given an exact solution $u$ and a numerical solution $\widehat{u}$ over $N$ spatial grid points $x_i$ at time $t_j$, the relative $L^1$ error is computed as
\begin{equation}
    \|u - \widehat{u} \|_1 = \frac{1}{N} \sum_{i = 1}^N \frac{| u(x_i, t_j) - \widehat{u}_i^j |}{| u(x_i, t_j) |}
\end{equation}


\subsubsection{Stability}

\begin{figure}
    \centering
    \begin{subfigure}[b]{0.48\textwidth}
        \centering
        \includegraphics[scale=0.45]{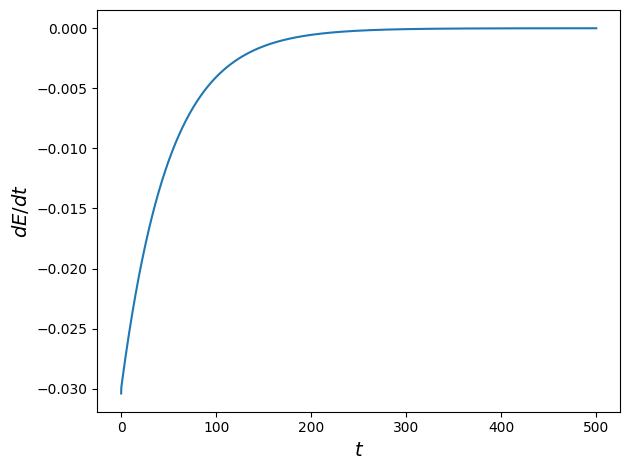}
        \caption{expansion model}
        \label{fig:stability}
    \end{subfigure}
    \begin{subfigure}[b]{0.48\textwidth}
        \centering
        \includegraphics[scale=0.45]{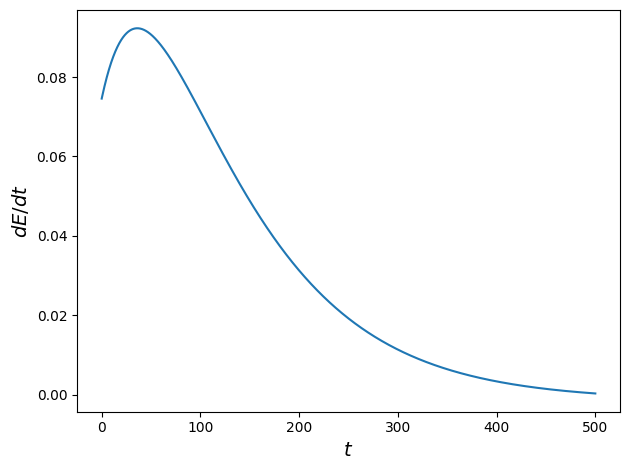}
        \caption{Fick-Jacobs model}
        \label{fig:fj-stability}
    \end{subfigure}
    \caption{Discrete $L^2$ energy derivative over time. (a) The energy of the expansion model system monotonically decays to an energy stable state, indicating stability over time. (b) The standard Fick-Jacobs model injects unphysical energy early on before $dE/dt \to 0$ from above, indicating it converges to a corrupted solution.}
\end{figure}

We verify the energy stability of the method. With $w(x) \equiv 1$, the scheme is the pure expansion-model discretization of Lemma~\ref{stab-lemma}. The test geometry $R(x) = 1 + \alpha \sin(kx)$ with $\alpha = 2$, $k = 1$ on the domain bounded by the first two positive roots of $\alpha\cos(2kx) - \sin(kx) = 0$ gives $R_{\min} \approx 2.19$, $\|R\|_\infty \approx 2.98$, and $\|R_x\|_\infty \approx 1.61$. The Lemma~\ref{stab-lemma} threshold requires $\Delta x \leq R_{\min}^2/(\|R\|_\infty^2 + \|R\|_\infty\|R_x\|_\infty) \approx 0.35$, comfortably satisfied by the grid spacings used here. We compute $dE/dt = \langle u, Lu\rangle_h + \langle u, f\rangle_h$ for both schemes. Figure~\ref{fig:stability} shows $dE/dt < 0$ for the expansion model throughout, asymptotically approaching zero from below as the solution relaxes. Figure~\ref{fig:fj-stability} shows the standard Fick-Jacobs discretization with $dE/dt > 0$ during an extended initial transient before $dE/dt \to 0$ from above, reflecting unphysical energy injection on a geometry with $\max|R_x| = 1.61$ outside the FJ range of validity \cite{Zwanzig1992}, leading the discretization to stabilize at a corrupted solution.

\subsubsection{Convergence}

\begin{figure}
    \centering
    \includegraphics[scale=0.5]{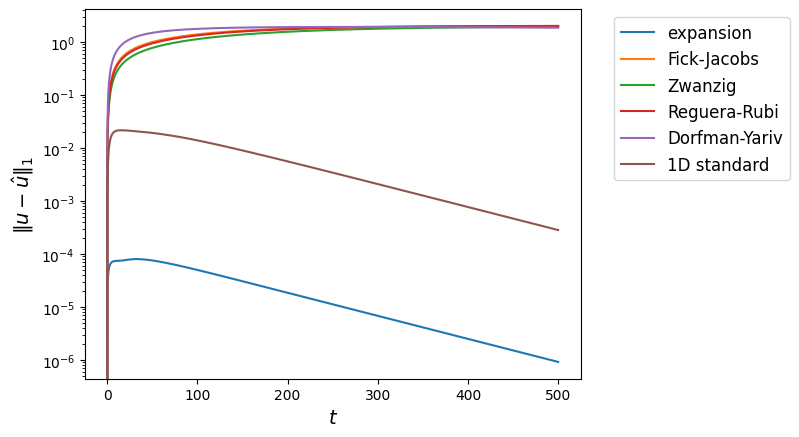}
    \caption{$L^1$ error over time for $N=80$, $\alpha = 2$. Only the expansion model is converges to the correct solution, while other models converge to an incorrect solution.}
    \label{fig:l1-conv}
\end{figure}

We now demonstrate the spatial convergence of our method. In Fig.~\ref{fig:l1-conv}, we see that our method is converging to the correct solution while the discretized Fick-Jacobs model and corrections are all converging to the same physically incorrect solution, confirming the conclusions from the previous section. As this solution should decay to a fixed constant concentration, the standard 1D model will likely eventually approach the same solution, which we also see in Fig.~\ref{fig:l1-conv}. The coefficient corrections produce models that are inconsistent with the true geometry-aware equation, and therefore cannot converge to the correct solution regardless of spatial resolution. This demonstrates that our model, unlike the corrections considered here, is consistent with the true geometry-aware equation and converges to the correct physical solution under grid refinement.

An important note is that there are several modifications are not included in Fig.~\ref{fig:l1-conv}, specifically under the Kalinay-Percus and Kalinay models. This is because the modifications to the diffusion coefficient in these cases require a division by $R_x$, which attains values numerically close to zero in the domain considered here. This meant that even the relatively modest problem was unable to be solved, demonstrating that these particular corrections are unable to complete simulations at the moderate radial gradients considered here.

\begin{figure}
    \centering
    \begin{subfigure}[b]{0.48\textwidth}
        \includegraphics[scale=0.45]{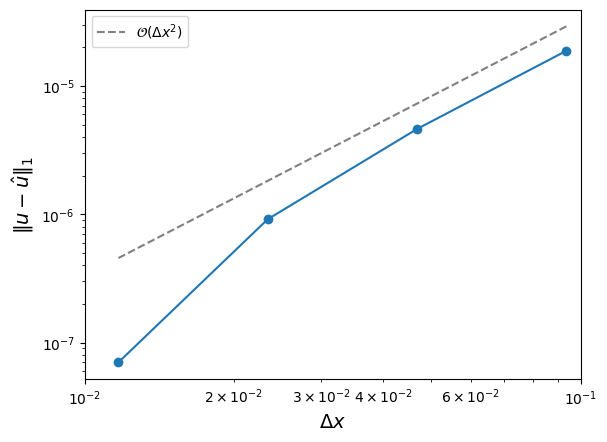}
        \caption{}
        \label{fig:dx-conv}
    \end{subfigure}
    \begin{subfigure}[b]{0.48\textwidth}
        \includegraphics[scale=0.45]{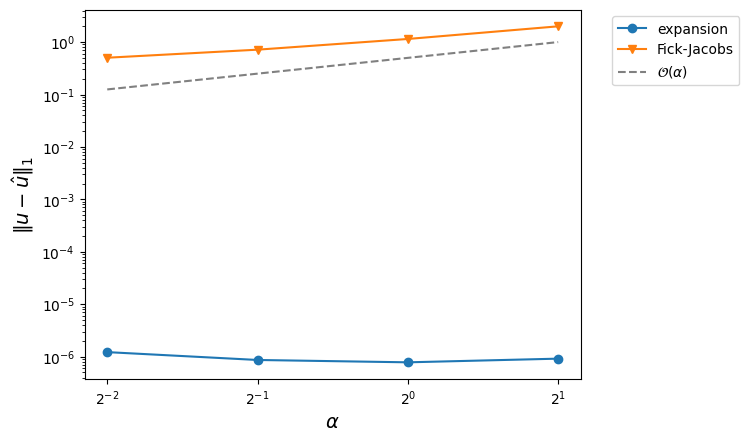}
        \caption{}
        \label{fig:alpha-conv}
    \end{subfigure}
    \caption{$L^1$ error convergence. Only the expansion model converges to the correct geometry-aware solution, while other models converge to an incorrect solution.}
\end{figure}

We next consider convergence rate. The test problem uses $w(x) \equiv 1$ and satisfies assumptions (A1)--(A2), so Theorem~\ref{thm:convergence} applies directly. In Fig.~\ref{fig:dx-conv}, we see that our method is converging at or above the expected model error rate of $\mathcal{O}(\Delta x^2)$, consistent with Theorem~\ref{thm:convergence}. We also see in Fig.~\ref{fig:alpha-conv} that the error remains relatively constant with changing $\alpha$. This is exactly the behavior that the derivation was targeted towards, and confirms that our model error is independent of $\alpha$. Within the triangle decomposition \eqref{eq:triangle-decomp}, term~II is $\mathcal{O}(\Delta x^2)$ uniformly in $\alpha$ by Theorem~\ref{thm:convergence}; the $\alpha$-independence of the total error reflects $\alpha$-independence of the modeling error (term~I), the defining property of the expansion model.

Together, these results show that the expansion model converges to the manufactured solution under refinement, while the standard diffusion-coefficient corrections converge to a different limit at the same $\Delta x \to 0$. The truncation in Jacobs' original derivation, and in particular the dropped volumetric-correction term on the left-hand side, is not recovered by any diffusion-coefficient modification of the form considered in the literature, so refinement of those discretizations cannot recover the geometry-aware limit, consistent with the structural analysis of Sec.~\ref{sec:original-fick-jacobs}. The expansion model retains the dropped terms by construction, and converges to the geometry-aware limit accordingly.


\subsection{Benchmark against axisymmetric ground truth on the truncated cone}\label{sec:bpb-benchmark}

The manufactured-solution test establishes numerical fidelity of the expansion model to its reduced PDE. We now turn to the complementary test of \emph{model fidelity to the underlying 3D physics} via a direct comparison against the three-dimensional axisymmetric diffusion equation.

The canonical benchmark is the truncated cone of Berezhkovskii, Pustovoit, and Bezrukov \cite{BPB2007} (hereafter BPB), defined as a long conical tube $R(x) = 1 + \lambda x$ on $[0, L]$ with $L = 20$, $D_0 = 1$, and constant wall slope $\lambda := R_x \in [0, 2]$. The cone isolates the dependence on $|R_x|$ alone ($R_{xx} = 0$). BPB ran Brownian dynamics simulations with reflecting walls and reflecting/absorbing end caps, and reported mean first-passage times $\tau_\lambda(n \to w)$ and $\tau_\lambda(w \to n)$ from each end cap to the other as functions of $\lambda$. We focus here on the $n \to w$ direction, the entropic-drift direction in which a particle moves from the narrow to the wide end of the cone. The opposite $w \to n$ direction, in which the entropic potential opposes net transport, exposes a limitation of the present form of $w(x)$ that we treat as an open problem and discuss in Sec.~\ref{sec:bpb-followup} below.

BPB established that the Reguera-Rub\'{i} correction accurately reproduces the simulation MFPT for $|R_x| \lesssim 1$ but diverges from the Brownian-dynamics ground truth for $|R_x| > 1$, while  Zwanzig and Fick-Jacobs perform worse. None of these corrections contains a free parameter that could extend its range of validity, and BPB concluded that the conventional one-dimensional reduction is reliable only when $|R_x| \lesssim 1$.

We compute the steady-state MFPT predicted by the expansion model in the $n \to w$ direction on the BPB cone geometry across $\lambda \in [0, 2]$ and into the large-slope regime $\lambda > 2$. Rather than solve the stationary backward equation directly, we use the equivalent survival-probability formulation, which fits naturally into the forward-in-time SBDF2 scheme of Sec.~\ref{sec:discretization} without modification. Specifically, we time-march \eqref{eq:dFJ} with the absorbing-cap boundary condition $u = 0$ at the wide end, the reflecting condition $F_{-1/2} = 0$ at the narrow end, and initial data concentrated as a unit mass on the reflecting (narrow) cap: $u(x_0, 0) = 1/\mathcal{V}_0$, $u(x_i, 0) = 0$ for $i > 0$, where $\mathcal{V}_0$ is the volume of the first cell. Letting $S(t) := \sum_i u_i(t)\, \mathcal{V}_i$ denote the surviving total mass (normalized so $S(0) = 1$), the cap-averaged MFPT is recovered via the standard identity $\tau_\lambda(n \to w) = \int_0^\infty S(t)\, dt$, which we evaluate numerically by trapezoidal integration on the recorded time history. We use $T = 2000$ with $\Delta t = 0.5$ throughout (4000 time steps), which is long enough at every $\lambda$ in the table for the survival probability to decay to numerical zero ($S(T) \lesssim 10^{-16}$); the tail contribution beyond $T$ is therefore negligible at any precision relevant to this comparison.

For the FJ and RR predictions reported in Table~\ref{tab:bpb-benchmark}, we do not solve the corresponding 1D PDEs numerically. The manufactured-solution test of Sec.~\ref{sec:numerical-results} has already established that the discretizations of FJ, Zwanzig, RR, and Dorfman-Yariv converge to physically incorrect solutions of the geometry-aware problem under refinement, so numerical MFPTs from those models would conflate two distinct sources of error and obscure the comparison. We therefore use the closed-form MFPT expression \cite{BPB2007} for any conventional 1D reduction with constant effective diffusivity $D_{\rm eff}$ on the cone $R(x) = 1 + \lambda x$,
\begin{equation}
    \tau_\lambda(n \to w) = \frac{L^2}{6 D_{\rm eff}} \cdot \frac{3 + \lambda L}{1 + \lambda L}, \label{eq:bpb-mfpt-formulas}
\end{equation}
evaluated with $D_{FJ}(x)$ and $D_{RR}(x)$ \cite{BPB2007}. These are the exact MFPTs that the FJ and RR equations analytically produce on the cone.

The two routes to the MFPT used here differ in numerical strategy. The BPB column solves the stationary backward MFPT equation directly via FEM on the 2D axisymmetric domain, which is the natural choice for a 2D elliptic problem where a single linear solve produces the full $\tau(x, r)$ field. The expansion column instead uses the equivalent survival-probability formulation $\tau = \int_0^\infty S(t)\, dt$ via the existing SBDF2 forward-in-time scheme of Sec.~\ref{sec:discretization}, which is the natural choice in our setting because the 1D scheme of \eqref{eq:dFJ} is already implemented as a forward time-marcher, and introducing a separate stationary-BVP solver would duplicate functionality. The two formulations are mathematically equivalent, since the MFPT is a property of the geometry and the PDE, not of the numerical method used to extract it, and the agreement between the columns at every $\lambda$ in Table~\ref{tab:bpb-benchmark} is a consequence of that equivalence rather than a coincidence of methods.

The $w_\delta$ values shown in the second column of Table~\ref{tab:bpb-benchmark} were obtained per $\lambda$ by manual tuning of the 1D solver to match the BPB ground truth in the $n \to w$ direction. A principled construction of $w(x)$ from the 3D physics that removes the need for per-geometry tuning, and the extension of the framework to handle the $w \to n$ direction in which entropic drift opposes net transport, are treated as open problems for follow-up work (see Sec.~\ref{sec:bpb-followup}).

\subsubsection{Results}

Table~\ref{tab:bpb-benchmark} reports the MFPT ratios $2D \tau_\lambda(n \to w)/L^2$ for the expansion model at per-$\lambda$ manually tuned $w_\delta$, alongside the BPB axisymmetric ground truth and the closed-form FJ and RR predictions from \eqref{eq:bpb-mfpt-formulas}.

Two observations follow. First, within the BPB range $\lambda \in [0, 2]$, the FJ closed-form MFPT underestimates the ground truth across the full range and the RR closed-form MFPT diverges progressively from the ground truth as $\lambda$ increases past unity, consistent with the BPB conclusion that the conventional 1D reductions are reliable only when $|R_x| \lesssim 1$. The tuned expansion model matches the ground truth at every $\lambda$ in the table, because the tuning is performed per-$\lambda$ to enforce that match, demonstrating that for each $\lambda$ a value of $w_\delta$ exists at which the expansion model reproduces the axisymmetric MFPT, while no value of any parameter can produce the same agreement for standard correction models.

Second, the wall-time columns $t_{\rm BPB}$ and $t_{\rm exp}$ quantify the cost of the two routes to the same MFPT. The BPB column requires an axisymmetric 3D FEM solve whose cost grows with the cone volume, from $\sim\!0.4$~s at $\lambda = 0$ to $\sim\!19$~s at $\lambda = 2$ (a factor of $\sim\!50$ across the range). The expansion-model column requires a 4000-step SBDF2 time-march of the 1D survival-probability problem, with constant cost $\sim\!1.3$~s independent of slope: the 1D grid does not grow with $\lambda$, the time horizon and step size are fixed, and only the geometry-dependent coefficients in \eqref{eq:dFJ} change between solves. The ratio $t_{\rm BPB}/t_{\rm exp}$ thus grows from $\sim\!0.3$ at $\lambda = 0$ (where the 3D problem is cheapest and the time-march has not yet broken even against a single FEM solve) to $\sim\!15$ at $\lambda = 2$. The reported wall times do not include the cost of tuning $w_\delta$; once a value is fixed for a given $\lambda$, it is reused without modification across all subsequent uses of the model on that geometry.

\subsubsection{Scope and open problems}\label{sec:bpb-followup}

Two limitations of the present benchmark are worth stating plainly. First, the tuning is performed per-$\lambda$ rather than once for the whole geometry family. A single $w_\delta$ for the cone --- and, more generally, a principled construction of $w(x)$ from the 3D physics that requires no per-geometry tuning --- is a substantial separate research program and is treated in follow-up work. Second, the benchmark above is restricted to the $n \to w$ direction. In the $w \to n$ direction, entropic drift opposes net transport rather than supporting it, and the present functional form of $w(x)$ (Eq.~\eqref{eq:dFJ}), which depends only on the magnitude $R_x^2$, does not distinguish between the two cases. The demonstration here establishes the foundational result that the expansion model has structural degrees of freedom which the no-parameter corrections lack, and which suffice to reproduce the axisymmetric ground truth in the entropic-drift direction across $|R_x| > 1$.

\begin{table}[h]
    \centering
    \begin{tabular}{c|c|cc|ccc|c}
        $\lambda$ & $w_\delta$ & BPB \cite{BPB2007} & $t_{\rm BPB}$ (s) & FJ & RR & expansion & $t_{\rm exp}$ (s) \\
        \hline
        0.00 & ---    & 1.000 & 0.41  & 1.000 & 1.000 & 1.000 & 1.30 \\
        0.25 & 0.0425 & 0.457 & 1.73  & 0.444 & 0.458 & 0.457 & 1.30 \\
        0.50 & 0.114  & 0.433 & 3.55  & 0.394 & 0.440 & 0.433 & 1.30 \\
        0.75 & 0.205  & 0.447 & 5.91  & 0.375 & 0.469 & 0.447 & 1.30 \\
        1.00 & 0.3095 & 0.473 & 7.84  & 0.365 & 0.516 & 0.473 & 1.30 \\
        1.25 & 0.423  & 0.502 & 10.07 & 0.359 & 0.575 & 0.502 & 1.30 \\
        1.50 & 0.541  & 0.530 & 13.63 & 0.355 & 0.640 & 0.530 & 1.30 \\
        1.75 & 0.6645 & 0.557 & 16.41 & 0.352 & 0.709 & 0.557 & 1.30 \\
        2.00 & 0.791  & 0.582 & 18.93 & 0.350 & 0.782 & 0.582 & 1.30 \\
    \end{tabular}
    \caption{Normalized $n \to w$ MFPT ratios $2 D \tau_\lambda(n \to w)/L^2$ on the BPB cone geometry ($R_0 = 1$, $L = 20$, $D_0 = 1$). The BPB column is an axisymmetric 3D ground truth obtained by reducing the 3D diffusion equation to a 2D axisymmetric finite-element MFPT problem; FJ and RR are evaluated in closed form via \eqref{eq:bpb-mfpt-formulas}; the expansion-model entries are obtained at the per-$\lambda$ manually tuned $w_\delta$ values shown in the second column. Wall times $t_{\rm BPB}$ report a single stationary MFPT solve (mesh generation, assembly, and LU solve) on the 2D axisymmetric FEM mesh; $t_{\rm exp}$ reports an SBDF2 time-march of the 1D expansion-model survival-probability problem with $T = 2000$ and $\Delta t = 0.5$ (4000 steps), plus initial-data setup and trapezoidal integration of $S(t)$. At $\lambda = 0$ the weight $w(x) = R_x^2/(R_x^2 + w_\delta^2)$ is identically zero regardless of $w_\delta$, so no tuned value is reported.}
    \label{tab:bpb-benchmark}
\end{table}


\subsection{Applying the F-J expansion model to a 3D biological system}

\begin{figure}
    \centering
    \includegraphics[scale=0.5]{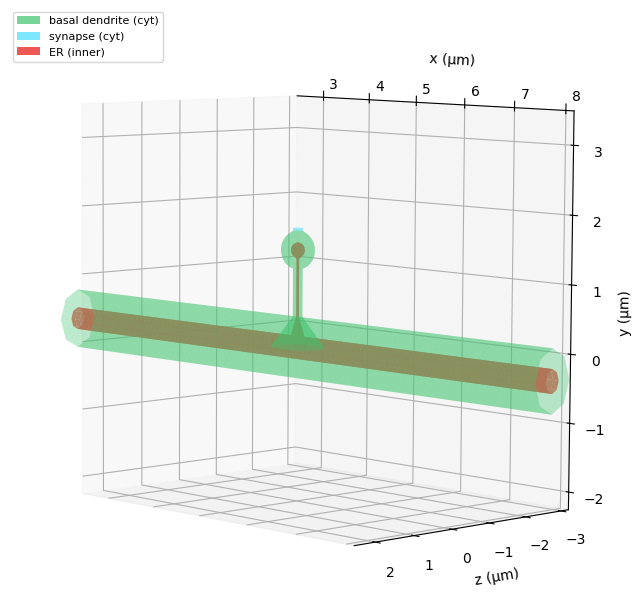}
    \caption{3D organization of a portion of a neuronal dendrite with a synaptic spine protrusion consisting of a spine neck and head (light green). Inside the neuron, the endoplasmic reticulum, a cell organelle, relevant for calcium exchange is visualized in red.}
    \label{fig:spine-geom}
\end{figure}

We now apply our expansion model to a highly relevant neurobiological system. When studying structure-function interplay in neurons, a neuron and its intracellular organelles need to be resolved in three space dimensions \cite{BG2021}. This is also true for contact regions between neurons, so called chemical synapses, that can be thought of as protrusions from the neuron's tube-like dendrites and axons \cite{Breit2016, Rosado2022}. The geometric organization of our model system is visualized in Fig.~\ref{fig:spine-geom}. Prior work demonstrated that standard 1D simulation is insufficient to capture spine dynamics and concluded that higher dimensional simulation is necessary \cite{MCENG2022, CS2017}. We demonstrate that this conclusion, while accurate for existing 1D models, does not hold for a properly derived geometry-aware reduction. 

In 1D, we represent the geometry in Fig.~\ref{fig:spine-geom} as a mesh consisting of two cables attached at a single branch point where each node contains geometric information about both the cytoplasm (light green) and endoplasmic reticulum organelle (red).

We then model the calcium dynamics, which play a crucial role in cell signaling, memory formation, learning, and aging. Calcium enters the spine head through a region called the postsynaptic density where, depending on the strength of the signal, it may travel through the head, down the neck, and into the dendrite. The full dynamics are modeled through a complex array of biological exchange mechanisms that act as surface fluxes on both the plasma membrane (outer domain boundary surface) and endoplasmic reticulum membrane (inner domain boundary surface). We do not describe the full biological model details here, but they may be found in \cite{BRNQ2023}.

For our test, we inject a 1 millisecond calcium stimulus and track the evolution of the calcium signal over 10 milliseconds by measuring the average cytosolic calcium concentration in the head, neck, and the section of the dendrite adjacent to the spine neck. We simultaneously run an equivalent 3D simulation with the same geometry and biological parameters using the uG4 framework \cite{VRRNW2013} to establish a ground truth, as well as compare to dynamics produced using the current standard 1D model.

\begin{figure}
    \centering
    \begin{subfigure}[b]{0.31\textwidth}
        \includegraphics[scale=0.3]{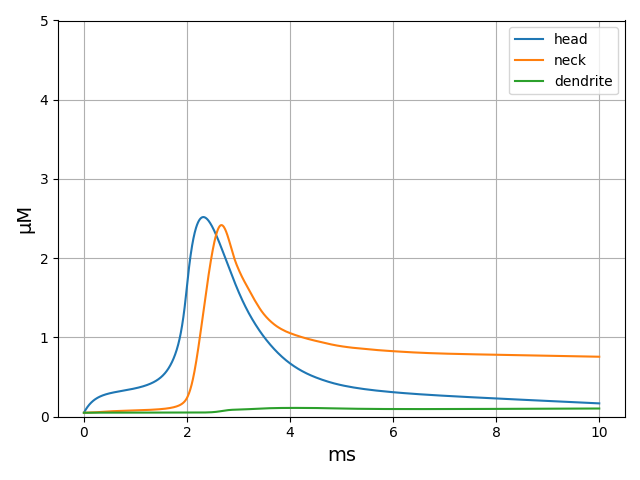}
        \caption{full 3D}
    \end{subfigure}
    \begin{subfigure}[b]{0.31\textwidth}
        \includegraphics[scale=0.3]{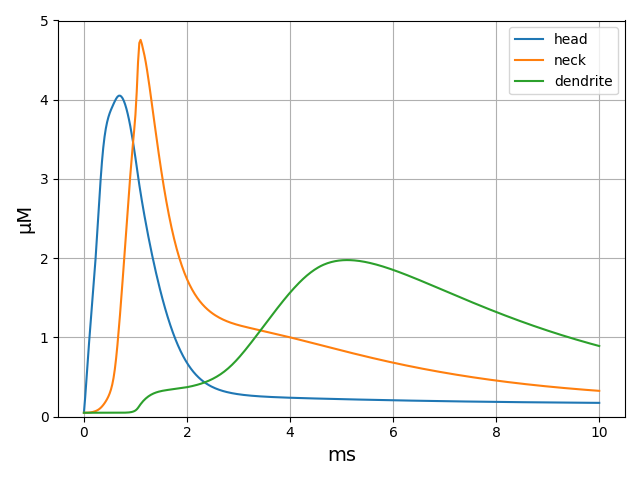}
        \caption{1D standard}
    \end{subfigure}
    \begin{subfigure}[b]{0.31\textwidth}
        \includegraphics[scale=0.3]{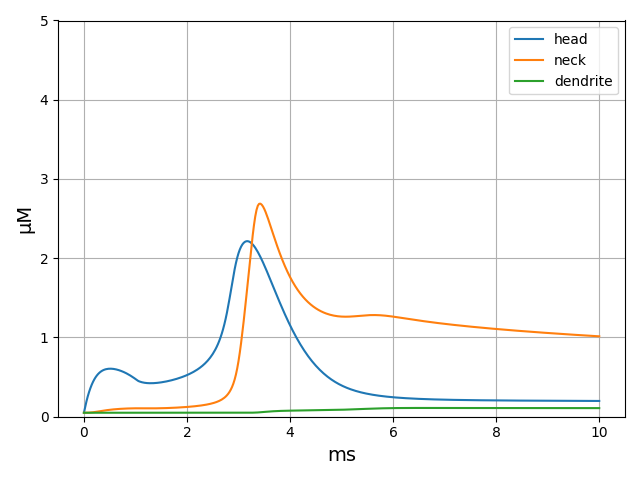}
        \caption{1D Fick-Jacobs expansion}
    \end{subfigure}
    \caption{Regional concentration of stimulated spine simulation using a full 3D simulation, the standard 1D diffusion reduction, and the expanded Fick-Jacobs model.}
    \label{fig:spine-data}
\end{figure}

We plot the measured concentrations in Fig.~\ref{fig:spine-data}. In the ground truth 3D simulation, a delayed small calcium spike is generated in the head and into the neck, but does not reach into the dendrite and no dendritic calcium wave is initiated. The standard 1D model fundamentally misrepresents these dynamics -- it generates a rapid, large-amplitude spike in the head that quickly travels into the neck and results in a calcium wave traveling through the dendrite. This is significant in that the spurious calcium wave in the dendrite can have far-reaching effects within a neuron.

With our model, we must initially use the ground truth results to tune the geometry-dependent parameter to an appropriate value, namely $w_{\delta, cyt} = 0.625$, $w_{\delta, er} = 0.85$. Once tuned, we see that our model has accurately reproduced the dynamics within the spine -- the delayed small spikes in the head and neck, and the suppression of dendritic wave propagation. This is significant in that a 1D model is able to faithfully produce physical dynamics that previously required a full (and costly) 3D simulation.

While our model required tuning of $w_\delta$ to the ground truth simulation, this is a one time cost that enables far-reaching biological and clinical applications. The 10 ms simulation took approximately 6 hours to complete using 8 processors on a desktop workstation, making biological parameter sweeps, clinical length simulations at the minute scale, and placement of many synapses on a complete 1D neuron model computationally infeasible without significant resources. In stark contrast, the 1D simulation took approximately 1 second to complete on a standard laptop. After the one time cost of parameter tuning, all of these simulations become reasonable in time scale, allowing for the kinds of extensive biological parameter studies, clinical timescale simulations, and multiscale neuronal modeling that the 3D approach renders computationally infeasible.



\section{Conclusion}

In this paper, we have returned to the foundations of Jacobs' original derivation and, through elementary analysis of the complete equation, identified two sources of geometric information loss — the truncation of higher order spatial terms and the loss of volumetric geometric information on the left hand side — that are inherited by all standard corrections from the literature. We have shown that these losses correspond to discarded $\Delta x^2$ terms in a systematic reduction of the 3D axisymmetric diffusion equation, and that retaining them through the next order in $\Delta x$ produces a geometry-aware model whose model error is independent of $R_x$. By treating the derivation as a Taylor expansion around the center of the subdomain and retaining this geometric information, we have produced a geometry-aware model with a provably stable and convergent discretization. This discretization converges under refinement to the geometry-aware reduction it is derived from, and an elementary structural argument shows that the standard diffusion-coefficient corrections converge to a different limit at the same $\Delta x \to 0$.

We have established discrete energy stability and convergence of the numerical method analytically, confirmed these results through numerical examples, validated the model against three-dimensional axisymmetric ground truth on the canonical truncated-cone benchmark, and applied the method to a branched network in the biological application. On the BPB cone, with $w_\delta$ tuned per $\lambda$, the expansion model reproduces the axisymmetric ground truth across the full slope range including $|R_x| > 1$, demonstrating that the structural degrees of freedom of the expansion model are sufficient to capture the geometric dependence that no parameter of Fick-Jacobs, Zwanzig, Reguera-Rub\'{i}, Kalinay-Percus, or Dorfman-Yariv can recover. We apply our method to realistic neuronal dendritic spines and demonstrate that a geometry-aware one-dimensional reduction can faithfully reproduce full three-dimensional spine simulation results that the standard reduction cannot achieve, at a fraction of the computational cost. This opens the door to biological parameter studies, clinical timescale simulations, and multiscale neuronal modeling that remain computationally intractable in three dimensions. Forthcoming work applies this framework to the study of calcium signaling and propagation in spines, building directly on the biological application presented here.

The geometry-aware framework developed here is broadly applicable to diffusion in tubular networks across biological, engineering, and physical applications. The blending function $w(x)$ enters as a modeling construct with the parameter $w_\delta$ currently set by tuning against ground-truth data. A principled construction of $w$ from the 3D physics --- without reference simulation --- is a substantial separate research program. Further development would have massive impact on simulating biological \cite{Grein2014,Breit2016,Breit2018,Rosado2022, Shirinpour2021,YKA1991,BRNQ2023} and other systems without the need for reference simulations.



\vspace{2em}
\noindent \textbf{Data availability.} Code and data for all examples are available at \url{https://github.com/zmiksis/miniSPINE}.\\
\noindent \textbf{Authors’ contributions.} Z.M.M.: conceptualization, formal analysis, investigation, software, validation, writing -- original draft, G.Q.: conceptualization, formal analysis, funding acquisition, investigation, methodology, project administration, supervision, writing -- review \& editing \\



\bibliographystyle{vancouver}
\bibliography{ref}

\end{document}